\documentclass[12pt]{article}
\usepackage{booktabs}
\usepackage{url}
\usepackage[figuresright]{rotating}

\usepackage{enumerate}

\RequirePackage[authoryear]{natbib}
\RequirePackage[colorlinks,citecolor=blue,urlcolor=blue]{hyperref}

\usepackage{bbm}
\usepackage{rotating}

\usepackage{amsmath}
\usepackage{multirow}
\usepackage{latexsym,amssymb,amsmath}
\usepackage{CJK}
\usepackage{picinpar,graphicx}
\usepackage{float}
\usepackage{listings}
\usepackage{color}
\usepackage{bm}
\usepackage{mathrsfs}
\usepackage{graphicx}
\usepackage{subfigure}
\usepackage{epstopdf}

\usepackage{multirow}

\usepackage{authblk}

\newtheorem{theorem}{Theorem}
\newtheorem{lemma}{Lemma}

\newtheorem{remark}{Remark}

\newcommand{\ba}{\begin{array}}
\newcommand{\ea}{\end{array}}
\newcommand{\bt}{\begin{tabular}}
\newcommand{\et}{\end{tabular}}
\newcommand{\btb}{\begin{table}}
\newcommand{\etb}{\end{table}}
\newcommand{\bc}{\begin{center}}
\newcommand{\ec}{\end{center}}
\newcommand{\bea}{\begin{eqnarray}}
\newcommand{\eea}{\end{eqnarray}}
\newcommand{\Bea}{\begin{eqnarray*}}
\newcommand{\Eea}{\end{eqnarray*}}
\newcommand{\beq}{\begin{equation}}
\newcommand{\eeq}{\end{equation}}
\newcommand{\ben}{\begin{enumerate}}
\newcommand{\een}{\end{enumerate}}

\def\T{{ \mathrm{\scriptscriptstyle T} }}

\begin{document}
\title{A studentized permutation test in group sequential designs}

\author[1]{Long-Hao Xu\thanks{Corresponding author: long-hao.xu@med.uni-goettingen.de}}
\author[2]{Tobias M\"{u}tze}
\author[3]{Frank Konietschke}
\author[1]{Tim Friede}

\affil[1]{Department of Medical Statistics, University Medical Center G\"{o}ttingen, G\"{o}ttingen, Germany}

\affil[2]{Statistical Methodology, Novartis Pharma AG, Basel, Switzerland}

\affil[3]{Institute of Biometry and Clinical Epidemiology, Charit\'{e} -- Universit\"{a}tsmedizin Berlin, Berlin, Germany}

\maketitle

\begin{abstract}
In group sequential designs, where several data looks are conducted for early stopping, we generally assume the vector of test statistics from the sequential analyses follows (at least approximately or asymptotially) a multivariate normal distribution. However, it is well-known that test statistics for which an asymptotic distribution is derived may suffer from poor small sample approximation. This might become even worse with an increasing number of data looks. The aim of this paper is to improve the small sample behaviour of group sequential designs while maintaining the same asymptotic properties as classical group sequential designs. This improvement is achieved through the application of a modified permutation test. In particular, this paper shows that the permutation distribution approximates the distribution of the test statistics not only under the null hypothesis but also under the alternative hypothesis, resulting in an asymptotically valid permutation test. An extensive simulation study shows that the proposed permutation test better controls the Type I error rate than its competitors in the case of small sample sizes.
\end{abstract}


\section{Introduction}

Group sequential designs are commonly used in clinical trials.
The main characteristic of this design is that the sample size is not fixed in advance and observations are sequentially obtained.
The benefits of group sequential design are widely acknowledged; a discussion can be found in \cite{JennisonTurnbull:1999,Proschan:2006,WassmerBrannath:2016,BhattMehta:2016,Neumannetal:2017}.
For large sample sizes, the classical methods of group sequential designs work well.
However, there are many situations in medical research where sample sizes are rather restricted.
For example, the sample size is usually small in preclinical experiments and early clinical trials such as first-in-human studies.
In addition, we often cannot obtain sufficiently large data sets in rare diseases.
From a statistical perspective, the above situations lead to small sample size problems.

The existing methods do not work well for small samples in that the Type I error rate may be inflated resulting in wrong conclusions; e.g., see Figure 1 in \cite{ShaoFeng:2007}.
There are some methods focusing on the small sample sizes problem in the case of the normal distribution.
For example, \cite{ShaoFeng:2007} computed critical values using Monte Carlo methods, whereas \cite{RomMcTague:2020} derived exact critical values through integration techniques.
In addition, \citet[Section 4.4]{JennisonTurnbull:1999}, \citet[Page 108]{WassmerBrannath:2016} and \cite{Nikolakopoulostal:2018} discussed a $t$-approximation that was proposed to address small sample problems.
However, these methods rely on the distributional assumption on observations which may not be satisfied in practice.
This motivates us to develop a general method that is independent of the distributional assumption, and the permutation test might be a good choice.

Permutation tests have been shown to control the Type I error rate accurately in various situations without distributional assumptions on observations; see for example \cite{SalmasoPesarin:2010} and the references therein.
In particular, following the ideas of \cite{Janssen:1997}, \cite{JanssenPauls:2003,JanssenPauls:2005}, \cite{NeubertBrunner:2007}, \cite{Pauly:2011}, \cite{KonietschkePauly:2012}, \cite{OmelkaPauly:2012}, \cite{ChungRomano:2013}, \cite{KonietschkePauly:2014}, \cite{Paulyetal:2015}, \cite{Paulyetal:2016}, and \cite{KonietschkeFriedePauly:2019}, we suggest a permutation test that is based on studentized statistics.


To date there are only few papers considering applications of permutation tests in group sequential designs.
For instance, \cite{Mehta:1994} only focused on linear rank statistics and mentioned that the corresponding large-sample theory is not very reliable for their method but did not provide any detailed investigations.
In addition, \citet[Section 8.2]{Proschan:2006} just briefly described the outline of permutation test in group sequential designs without providing further explanations or a theoretical basis.
Particularly, this method, as it is not based on studentized statistics, yielded non-robust results.
Although \citet[Section 6.2]{Mutzetal:2019} mentioned the use of permutation tests in group sequential designs, the correlation among these test statistics computed at the end of each stage was not considered.
Therefore, we develop a general robust method that applies the studentized permutation test for both small sample sizes and large sample sizes, which controls the Type I error rate well.
For this proposed method, we use the permutation distribution to approximate the distribution of the test statistics under the null hypothesis, and we compute the critical values based on the permutation distribution.

In this paper, we consider $K$-stage unbalanced designs for testing equality of the means between two arms.
We will carry out simulations to investigate finite sample size properties including Type I error rate and power, in particular for small samples.
For the asymptotic results, we will provide mathematical proofs showing that the conditional permutation distribution given the observed data approximates the unconditional distribution of the test statistics.
This implies that our proposed method is applicable to both small sample sizes and large sample sizes.
In the special case of one stage $K=1$, \citet{ChungRomano:2013} have stated a similar result for our proposed method.

The paper is organized as follows.
In Section \ref{sec-2} we review the classical group sequential designs and propose the test decision function of group sequential designs.
The novel permutated version of group sequential designs and its theoretical properties are stated and explained in Section \ref{sec-3}.
In Section \ref{sec-4}, both small sample behaviour and large sample behaviour of this permutation test are examined in an extensive simulation study.
Concluding remarks are presented in Section \ref{sec-conclusion}.
All proofs are given in the supplementary material.

\section{Classical group sequential designs}\label{sec-2}

Consider the two-arm comparison in the general $K$-stage unbalanced design.
Let $m_k$ be the number of the new observations for the treatment arm at the end of stage $k=1,\ldots,K$.
Similarly, let $n_k$ be the number of the new observations for the control arm at the end of stage $k=1,\ldots,K$.
At the end of stage $k$ we have $\widetilde{m}_k$ observations available for analysis in the treatment arm and $\widetilde{n}_k$ observations available for analysis in the control arm, where $\widetilde{m}_k=\sum_{j=1}^{k}m_j$ and $\widetilde{n}_k=\sum_{j=1}^{k}n_j$.
From a practical view, we assume the allocation ratio between the two arms is fixed, such that all $k=1,\ldots,K$,
\begin{align}\label{gamma}
\frac{m_k}{n_k}=\gamma>0.
\end{align}

Assume that $X_{i}\overset{iid}{\sim}F_1$ for $i=1,\ldots,\widetilde{m}_K$, and $Y_{j}\overset{iid}{\sim}F_2$ for $j=1,\ldots,\widetilde{n}_K$, are from the treatment arm and the control arm, respectively.
The true mean in each arm is denoted as $\mu_1=E(X_1)<\infty$ and $\mu_2=E(Y_1)<\infty$.
Similarly, the true variance in each arm is denoted as $\sigma_1^2=\mbox{var}(X_1)<\infty$ and $\sigma_2^2=\mbox{var}(Y_1)<\infty$.
To illustrate the characteristics of group sequential designs, we introduce the symbol $\biguplus$ to represent the concatenation of two data sets.
Specifically, we define $\biguplus_{j=1}^i{Z}_j=(X_{1},\ldots,X_{\widetilde{m}_i},Y_{1},\ldots,Y_{\widetilde{n}_i})$ for $i=1,\ldots,K$, where ${Z}_j$ represents the new observations at the end of stage $j=1,\ldots,i$, namely
\begin{align}\label{Z}
{Z}_j&=\left(X_{\widetilde{m}_{j-1}+1},\ldots,X_{\widetilde{m}_j},Y_{\widetilde{n}_{j-1}+1},\ldots,Y_{\widetilde{n}_j}\right)
\end{align}
with the definitions of $\widetilde{m}_{0}=0$ and $\widetilde{n}_{0}=0$.

Without loss of generality, we first consider testing the null hypothesis $H_0: \mu_1=\mu_2$ against the alternative hypothesis $H_1: \mu_1>\mu_2$.
Under the null hypothesis $H_0: \mu_1=\mu_2$, the corresponding studentized test statistic at the end of stage $k=1,\ldots,K$ is
\begin{align}\label{eq-statistics}
S_{\widetilde{m}_k,\widetilde{n}_k}
=S_{\widetilde{m}_k,\widetilde{n}_k}\left(\biguplus_{i=1}^k{Z}_i\right)
=\frac{\left(\frac{\widetilde{m}_k\widetilde{n}_k}{\widetilde{m}_k+\widetilde{n}_k}\right)^{1/2}\left(\widehat{\mu}_{1,k}-\widehat{\mu}_{2,k}\right)}{\left(\frac{\widetilde{n}_k}{\widetilde{m}_k+\widetilde{n}_k}\widehat{\sigma}^2_{1,k}+\frac{\widetilde{m}_k}{\widetilde{m}_k+\widetilde{n}_k}\widehat{\sigma}^2_{2,k}\right)^{1/2}},
\end{align}
where
\begin{align*}
\widehat{\mu}_{1,k} &= \widehat{\mu}_{1,k}\left(\biguplus_{i=1}^k{Z}_i\right)=\frac{1}{\widetilde{m}_k}\sum_{i=1}^{\widetilde{m}_k}X_i, \\
\widehat{\mu}_{2,k} &= \widehat{\mu}_{2,k}\left(\biguplus_{i=1}^k{Z}_i\right)=\frac{1}{\widetilde{n}_k}\sum_{i=1}^{\widetilde{n}_k}Y_i,\\
\widehat{\sigma}^2_{1,k} &= \widehat{\sigma}^2_{1,k}\left(\biguplus_{i=1}^k{Z}_i\right)=\frac{1}{\widetilde{m}_k-1}\sum_{i=1}^{\widetilde{m}_k}\left(X_i-\frac{1}{\widetilde{m}_k}\sum_{j=1}^{\widetilde{m}_k}X_j\right)^2,\\
\widehat{\sigma}^2_{2,k} &= \widehat{\sigma}^2_{2,k}\left(\biguplus_{i=1}^k{Z}_i\right)=\frac{1}{\widetilde{n}_k-1}\sum_{i=1}^{\widetilde{n}_k}\left(Y_i-\frac{1}{\widetilde{n}_k}\sum_{j=1}^{\widetilde{n}_k}Y_j\right)^2.
\end{align*}


In this sequential setting we utilize the available observations $X_{1},\ldots,X_{\widetilde{m}_k},Y_{1},\ldots,Y_{\widetilde{n}_k}$ to compute the test statistic $S_{\widetilde{m}_k,\widetilde{n}_k}$ at the end of stage $k=1,\ldots,K$, resulting in the sequence $\{S_{\widetilde{m}_1,\widetilde{n}_1},\ldots,S_{\widetilde{m}_K,\widetilde{n}_K}\}$.
Given the critical values $\{c_1,\ldots,c_K\}$, a one-sided hypothesis testing $H_0: \mu_1=\mu_2$ against $H_1: \mu_1>\mu_2$ for classical group sequential designs is described as follows. After stage $k=1,\ldots,K-1$, if $S_{\widetilde{m}_k,\widetilde{n}_k}\geq c_k$, this procedure stops, and $H_0$ is rejected; otherwise, this procedure continues to stage $k+1$. After stage $K$, if $S_{\widetilde{m}_K,\widetilde{n}_K}\geq c_K$, this procedure stops, and $H_0$ is rejected; otherwise, this procedure stops, and $H_0$ is not rejected.



Similar to the classical theory of group sequential designs (cf. \cite{Scharfstein:1997}; \cite{JennisonTurnbull:1997}; \cite{JennisonTurnbull:1999}), letting $n$ approach infinity means that an increasing number of individuals arrive in that fixed time frame. We assume that for each $k=1,\ldots,K$,
\begin{align}\label{assumption-GSD}
\lim_{n\rightarrow\infty}\widehat{\mathcal{I}}_k/n=\lim_{n\rightarrow\infty}\left\{n\widehat{\mbox{var}}(\widehat{\mu}_{1,k}-\widehat{\mu}_{2,k})\right\}^{-1}=\lim_{n\rightarrow\infty}\left\{n\mbox{var}(\widehat{\mu}_{1,k}-\widehat{\mu}_{2,k})\right\}^{-1}=\mathcal{I}_k^*>0,
\end{align}
where $\widehat{\mathcal{I}}_k=(\widehat{\mbox{var}}(\widehat{\mu}_{1,k}-\widehat{\mu}_{2,k}))^{-1}$ is the information accrued at the corresponding analysis $k$, $\widehat{\mbox{var}}(\widehat{\mu}_{1,k}-\widehat{\mu}_{2,k})$ denotes the estimated variance of $\widehat{\mu}_{1,k}-\widehat{\mu}_{2,k}$, $\mbox{var}(\widehat{\mu}_{1,k}-\widehat{\mu}_{2,k})$ denotes the true variance of $\widehat{\mu}_{1,k}-\widehat{\mu}_{2,k}$, and $\mathcal{I}_k^*$ is called the inverse of the asymptotic variance of $\widehat{\mu}_{1,k}-\widehat{\mu}_{2,k}$ (as $n\rightarrow\infty$).
Under the null hypothesis $H_0: \mu_1=\mu_2$, the vector $S=(S_{\widetilde{m}_1,\widetilde{n}_1},\ldots,S_{\widetilde{m}_K,\widetilde{n}_K})^\T$ asymptotically follows the multivariate normal distribution $N({0}_K,{\Sigma})$ where ${0}_K$ denotes the $K$-dimensional vector of zeros and the $(i,j)$th element of covariance matrix $\Sigma$ is
\begin{align*}
\Sigma_{ij}=\mbox{cov}\left(S_{\widetilde{m}_i,\widetilde{n}_i},S_{\widetilde{m}_j,\widetilde{n}_j}\right)=\left(\mathcal{I}_{\min\{i,j\}}^*/\mathcal{I}_{\max\{i,j\}}^*\right)^{1/2}
\end{align*}
for $i=1,\ldots,K$ and $j=1,\ldots,K$.


In classical group sequential designs, \cite{Pocock:1977} and \cite{OBrienFleming:1979} proposed two different types of critical values based on the above multivariate normal distribution $N({0}_K,{\Sigma})$.
Without a prespecified maximum number of stages $K$, a flexible method for choosing critical values involves the use of $\alpha$-spending function (cf. \cite{LanDeMets:1983}; \cite{KimDeMets:1987}; \cite{Hwang:1990}).
This means that a specific error rate is allocated to each stage according to an $\alpha$-spending function $f(\cdot)$.
It is important to note that the $\alpha$-spending function $f(t),\ 0\leq t\leq1$ can be any non-decreasing function with $f(0)=0$ and $f(1)=\alpha$.
In this paper, we will use Pocock type $\alpha$-spending function and O’Brien-Fleming type $\alpha$-spending function.
For the Pocock type $\alpha$-spending function, we take $f(t)=\min\{\alpha\log(1+(e-1)t),\alpha\}$.
For the O’Brien-Fleming type $\alpha$-spending function, we take $f(t)=\min\{2-2\Phi(\Phi^{-1}(1-\alpha/2)/\surd{t}),\alpha\}$ if it is a one-sided hypothesis testing, and we take $f(t)=\min\{4-4\Phi(\Phi^{-1}(1-\alpha/4)/\surd{t}),\alpha\}$ if it is a two-sided hypothesis testing, where $\Phi(\cdot)$ is the cumulative distribution function of the standard normal distribution and $\Phi^{-1}(\cdot)$ is its inverse function.

After choosing a suitable $\alpha$-spending function $f$, we prespecify the maximum information $\mathcal{I}_{\max}$ (defined in \citet[Section 4]{Scharfstein:1997} and \citet[Section 4]{VanLancker:2022}) and we find critical values $c_k$ that solves the equations
\begin{equation}
\label{c-i}
\left\{
\begin{aligned}
& \mbox{pr}\left(S_{\widetilde{m}_1,\widetilde{n}_1}\geq c_1\right) = f\left(\widehat{\mathcal{I}}_1/\mathcal{I}_{\max}\right) \\
& \mbox{pr}\left(S_{\widetilde{m}_1,\widetilde{n}_1}<c_1,S_{\widetilde{m}_2,\widetilde{n}_2}\geq c_2\right) = f\left(\widehat{\mathcal{I}}_2/\mathcal{I}_{\max}\right) - f\left(\widehat{\mathcal{I}}_{1}/\mathcal{I}_{\max}\right) \\
& \ldots \\
& \mbox{pr}\left(S_{\widetilde{m}_1,\widetilde{n}_1}<c_1,\ldots,S_{\widetilde{m}_{K-1},\widetilde{n}_{K-1}}<c_{K-1},S_{\widetilde{m}_K,\widetilde{n}_K}\geq c_K\right) = \alpha - f\left(\widehat{\mathcal{I}}_{K-1}/\mathcal{I}_{\max}\right),
\end{aligned}
\right.
\end{equation}
where the maximum information $\mathcal{I}_{\max}$ is usually determined such that the trial has a specific power (e.g., 90\%) for a certain effect size.
For design purposes, we often assume that
\begin{align*}
\mathcal{I}_k = \frac{k}{K}\mathcal{I}_{\max}, \ k=1,\ldots,K,
\end{align*}
which means that information levels will be equally spaced up to $\mathcal{I}_{\max}$.
For more details about classical group sequential designs, the readers are referred to \cite{JennisonTurnbull:1999,WassmerBrannath:2016}.

In hypothesis testing, the test decision function, also known as the critical function, is an important tool for analyzing Type I error rate and power (cf. \citet[Section 2.3]{Hajek:1999}; \cite{LehmannRomano:2005}).
To the best of our knowledge, the test decision function of group sequential designs has not been explicitly introduced in the literature.
Therefore, we offer a novel definition based on the characteristics of group sequential designs.
Define the test decision function of group sequential designs as
\begin{align}\label{phi}
\varphi_n = \sum_{i=1}^K\varphi_{n,i},
\end{align}
where
\begin{equation*}
\varphi_{n,1}= \left\{
    \begin{aligned}
    1 & \quad \mbox{if} \ S_{\widetilde{m}_1,\widetilde{n}_1}\geq c_1  \\
    0 & \quad \mbox{otherwise}
    \end{aligned}
    \right.
\end{equation*}
and
\begin{equation*}
\varphi_{n,i}= \left\{
    \begin{aligned}
    1 & \quad \mbox{if} \ S_{\widetilde{m}_1,\widetilde{n}_1}< c_1 \ \mbox{and} \ \ldots \ \mbox{and} \ S_{\widetilde{m}_{i-1},\widetilde{n}_{i-1}}< c_{i-1} \ \mbox{and} \ S_{\widetilde{m}_i,\widetilde{n}_i}\geq c_i  \\
    0 & \quad \mbox{otherwise}
    \end{aligned}
    \right.
\end{equation*}
for $i=2,\ldots,K$.
The critical values $\{c_1,\ldots,c_K\}$ in $\varphi_n$ are defined in (\ref{c-i}).
Since observations are sequentially obtained and interim analyses are conducted, we define $\varphi_n$ as the summation of each analysis $\varphi_{n,i}$.
With the definition of $\varphi_n$, it is straightforward to show ${E}(\varphi_n)\rightarrow\alpha$ under the null hypothesis and ${E}(\varphi_n)\rightarrow1$ under the alternative hypothesis as $n\rightarrow\infty$.

To address the small sample problem for the test statistics defined in (\ref{eq-statistics}), \cite{Welch:1938,Welch:1947} proposed a $t$-approximation to compute critical values in the non-sequential setting. As suggested by \citet[Section 4.4]{JennisonTurnbull:1999} and \citet[Page 108]{WassmerBrannath:2016}, in group sequential designs, we employ the same stage levels for the $t$-approximation and recalculate critical values as $\tilde{c}_k=G^{-1}_{df_k}(\Phi(c_k))$, where $G^{-1}_{df_k}$ denotes the inverse of the cumulative distribution function of the $t$-distribution with $df_k$ degrees of freedom and
\begin{align*}
df_k=\frac{ \left( \widehat{\sigma}^2_{1,k}/\widetilde{m}_k + \widehat{\sigma}^2_{2,k}/\widetilde{n}_k \right)^2 }{ \widehat{\sigma}^4_{1,k}/\{\widetilde{m}_k^2(\widetilde{m}_k-1)\} + \widehat{\sigma}^4_{2,k}/\{\widetilde{n}_k^2(\widetilde{n}_k-1)\} }, \ k=1,\ldots,K.
\end{align*}
The corresponding test decision function of this $t$-approximation method is defined as
\begin{align*}
\tilde{\varphi}_n = \sum_{i=1}^K\tilde{\varphi}_{n,i},
\end{align*}
where
\begin{equation*}
\tilde{\varphi}_{n,1}= \left\{
    \begin{aligned}
    1 & \quad \mbox{if} \ S_{\widetilde{m}_1,\widetilde{n}_1}\geq \tilde{c}_1  \\
    0 & \quad \mbox{otherwise}
    \end{aligned}
    \right.
\end{equation*}
and
\begin{equation*}
\tilde{\varphi}_{n,i}= \left\{
    \begin{aligned}
    1 & \quad \mbox{if} \ S_{\widetilde{m}_1,\widetilde{n}_1}< \tilde{c}_1 \ \mbox{and} \ \ldots \ \mbox{and} \ S_{\widetilde{m}_{i-1},\widetilde{n}_{i-1}}< \tilde{c}_{i-1} \ \mbox{and} \ S_{\widetilde{m}_i,\widetilde{n}_i}\geq \tilde{c}_i  \\
    0 & \quad \mbox{otherwise}
    \end{aligned}
    \right.
\end{equation*}
for $i=2,\ldots,K$.

In the following sections, we will explore the counterpart of $\varphi_{n}$ in the permutated version of group sequential designs and compare their properties.

\section{The proposed studentized permutation test}\label{sec-3}

In this section, we introduce the studentized permutation test for group sequential designs.
We will begin by outlining the general framework of this proposed method and subsequently show its asymptotic properties.

When applying the permutation test to group sequential designs, the challenge of permuting the observed data in this sequential setting arises.
In classical group sequential designs, observations are sequentially collected, and there exists covariance between $S_{\widetilde{m}_i,\widetilde{n}_i}$ and $S_{\widetilde{m}_j,\widetilde{n}_j}$ for $i=1,\ldots,K$ and $j=1,\ldots,K$.
However, the classical permutation test cannot be applied directly because it does not consider the sequential nature of group sequential designs. For this reason, we propose stage-wise permutation for permutating the observed data, inspired by \citet[Section 8.2]{Proschan:2006} and \cite{Proschan:2009}.
The main idea of this proposed stage-wise permutation is that the observations stay within their original stage after the permutation. Later, we will show that this procedure preserves the correlation structure between these test statistics computed from different stages.

With the definition of ${Z}_j$ in (\ref{Z}), let
\begin{align*}
{Z}_j &=\left(X_{\widetilde{m}_{j-1}+1},\ldots,X_{\widetilde{m}_j},Y_{\widetilde{n}_{j-1}+1},\ldots,Y_{\widetilde{n}_j}\right)=\left(Z_{j,1},\ldots,Z_{j,m_j+n_j}\right), \ j=1,\ldots,K.
\end{align*}
We adopt the notations used in \cite{Paulyetal:2015,Paulyetal:2016}.
Let $\tau_i$ be a random variable that is uniformly distributed on the symmetric group $\mathcal{S}_{N_i}$, that is, the set of all permutation of the numbers $1,\ldots,N_i$, being independent from the data ${Z}_i=(Z_1,\ldots,Z_{N_i})$, where $N_i=m_i+n_i$ for $i=1,\ldots,K$.
Assume all $\tau_i$ are independent.
Furthermore, we assume that the variables $X_i$ and $Y_j$ are defined on the probability space $(\Omega,\mathcal{A},\mbox{pr})$ and are independent from the random permutations $\tau_i$ which are defined on the probability space $({\Omega_1},{\mathcal{A}_1},\widetilde{\mbox{pr}})$.
Here independence means independence on the joint probability space $(\Omega\times{\Omega_1},\mathcal{A}\otimes{\mathcal{A}_1},\mbox{pr}\otimes\widetilde{\mbox{pr}})$, where all random variables can be defined in an obvious manner via projections.


Now, we introduce the method of stage-wise permutation in group sequential designs. Given the new observations ${Z}_k=(Z_{k,1},\ldots,Z_{k,N_k})$ for $k=1,\ldots,K$, we can represent the permuted samples as ${Z}_k^{\tau_k}=(Z_{k,\tau_k(1)},\ldots,Z_{k,\tau_k(N_k)})$, where the permutation is applied to the new observations at each stage. Then the resulting permuted samples ${Z}_k^{\tau_{k}}$ are combined with the results $\biguplus_{i=1}^{k-1}{Z}_i^{\tau_{i}}$ from the preceding $k-1$ stages to yield the results $\biguplus_{i=1}^k{Z}_i^{\tau_{i}}$ after completing $k$ stages.
Based on the above notations, the test statistic at the end of stage $k=1,\ldots,K$ is
\begin{align}\label{eq-statistics-tau}
S_{\widetilde{m}_k,\widetilde{n}_k}^{\tau}
=S_{\widetilde{m}_k,\widetilde{n}_k}\left(\biguplus_{i=1}^k{Z}_i^{\tau_{i}}\right)
=\frac{\left(\frac{\widetilde{m}_k\widetilde{n}_k}{\widetilde{m}_k+\widetilde{n}_k}\right)^{1/2}\left(\widehat{\mu}_{1,k}^\tau-\widehat{\mu}_{2,k}^\tau\right)}{\left(\frac{\widetilde{n}_k}{\widetilde{m}_k+\widetilde{n}_k}\widehat{\sigma}^2_{1,k,\tau}+\frac{\widetilde{m}_k}{\widetilde{m}_k+\widetilde{n}_k}\widehat{\sigma}^2_{2,k,\tau}\right)^{1/2}},
\end{align}
where
\begin{align*}
\widehat{\mu}_{1,k}^\tau &= \widehat{\mu}_{1,k}\left(\biguplus_{i=1}^k{Z}_i^{\tau_{i}}\right)=\frac{1}{\widetilde{m}_k}\sum_{s=1}^k\sum_{i=1}^{m_s}Z_{s,\tau_s(i)},\\
\widehat{\mu}_{2,k}^\tau &= \widehat{\mu}_{2,k}\left(\biguplus_{i=1}^k{Z}_i^{\tau_{i}}\right)=\frac{1}{\widetilde{n}_k}\sum_{s=1}^k\sum_{i=m_s+1}^{m_s+n_s}Z_{s,\tau_s(i)},\\
\widehat{\sigma}^2_{1,k,\tau} &= \widehat{\sigma}^2_{1,k}\left(\biguplus_{i=1}^k{Z}_i^{\tau_{i}}\right)=\frac{1}{\widetilde{m}_k-1}\sum_{u=1}^{k}\sum_{i=1}^{m_u}\left(Z_{u,\tau_u(i)}-\frac{1}{\widetilde{m}_k}\sum_{s=1}^{k}\sum_{j=1}^{m_s}Z_{s,\tau_s(j)}\right)^2,\\
\widehat{\sigma}^2_{2,k,\tau} &= \widehat{\sigma}^2_{2,k}\left(\biguplus_{i=1}^k{Z}_i^{\tau_{i}}\right)=\frac{1}{\widetilde{n}_k-1}\sum_{u=1}^{k}\sum_{i=m_u+1}^{m_u+n_u}\left(Z_{u,\tau_u(i)}-\frac{1}{\widetilde{n}_k}\sum_{s=1}^{k}\sum_{j=m_s+1}^{m_s+n_s}Z_{s,\tau_s(j)}\right)^2.
\end{align*}

These values $\widehat{\mu}_{1,k}^\tau$, $\widehat{\mu}_{2,k}^\tau$, $\widehat{\sigma}^2_{1,k,\tau}$, and $\widehat{\sigma}^2_{2,k,\tau}$ denote the quantities $\widehat{\mu}_{1,k}$, $\widehat{\mu}_{2,k}$, $\widehat{\sigma}^2_{1,k}$, and $\widehat{\sigma}^2_{2,k}$ being computed with the permuted variables.
Essentially, the test statistic defined in (\ref{eq-statistics-tau}) is computed from the permutated sample $\biguplus_{i=1}^k{Z}_i^{\tau_{i}}$ instead of the original sample $\biguplus_{i=1}^k{Z}_i$.

Next, we consider the computation of the corresponding critical values $c_k^\tau$ based on the test statistics $S_{\widetilde{m}_k,\widetilde{n}_k}^{\tau}$.
Similar to (\ref{c-i}), we adopt the $\alpha$-spending function approach and find critical values $c_k^\tau$ that solves the equations
\begin{equation}
\label{c-i-tau}
\left\{
\begin{aligned}
& \widetilde{\mbox{pr}}\left(S_{\widetilde{m}_1,\widetilde{n}_1}^{\tau}\geq c_1^\tau\right) = f\left(\widehat{\mathcal{I}}_1/\mathcal{I}_{\max}\right) \\
& \widetilde{\mbox{pr}}\left(S_{\widetilde{m}_1,\widetilde{n}_1}^{\tau}< c_1^\tau, S_{\widetilde{m}_2,\widetilde{n}_2}^{\tau}\geq c_2^\tau\right)  = f\left(\widehat{\mathcal{I}}_2/\mathcal{I}_{\max}\right) - f\left(\widehat{\mathcal{I}}_{1}/\mathcal{I}_{\max}\right) \\
& \ldots \\
& \widetilde{\mbox{pr}}\left(S_{\widetilde{m}_1,\widetilde{n}_1}^{\tau}< c_1^\tau, \ldots, S_{\widetilde{m}_{K-1},\widetilde{n}_{K-1}}^{\tau}< c_{K-1}^\tau, S_{\widetilde{m}_K,\widetilde{n}_K}^{\tau}\geq c_K^\tau\right) = \alpha - f\left(\widehat{\mathcal{I}}_{K-1}/\mathcal{I}_{\max}\right).
\end{aligned}
\right.
\end{equation}

In (\ref{c-i-tau}), we use the permutation distribution to approximate the distribution of the test statistics under the null hypothesis.
This permutation distribution is data-driven and reflects the characteristic of the observed data, particularly for small samples size.
Consequently, our proposed method provides a new way to generating critical values that depend on the observed data.

In the next steps, we will show the asymptotic properties of the proposed method.

\begin{theorem}\label{theorem-1}
Under the assumption of (\ref{assumption-GSD}), the conditional permutation distribution of
\begin{align*}
{S}^\tau=\left(S^{\tau}_{\widetilde{m}_1,\widetilde{n}_1},S^{\tau}_{\widetilde{m}_2,\widetilde{n}_2},\ldots,S^{\tau}_{\widetilde{m}_K,\widetilde{n}_K}\right)^\T
\end{align*}
given the observed data $\biguplus_{i=1}^K{Z}_i$ weakly converges to a multivariate normal distribution $N({0}_K,{\Sigma})$ in probability, where ${0}_K$ denotes the $K$-dimensional vector of zeros and the $(i,j)$th element of ${\Sigma}$ is
\begin{align*}
{\Sigma}_{ij}=\mbox{cov}\left(S^{\tau}_{\widetilde{m}_i,\widetilde{n}_i},S^{\tau}_{\widetilde{m}_j,\widetilde{n}_j}\right)=\left(\mathcal{I}^*_{\min\{i,j\}}/\mathcal{I}^*_{\max\{i,j\}}\right)^{1/2}
\end{align*}
for $i=1,\ldots,K$ and $j=1,\ldots,K$.
\end{theorem}

Similar to \cite{Paulyetal:2015}, the above convergence in terms of distances means
\begin{align}\label{s-tau-convergence}
\rho_K\left(\mathcal{L}\left({S}^\tau|\biguplus_{i=1}^K{Z}_i\right),N({0}_K,{\Sigma})\right)\overset{p}{\rightarrow}0,
\end{align}
where $\rho_K$ denotes a distance that metrizes weak convergence on the $K$-dimensional space, e.g. the Prohorov distance, see Section 11.3 in \cite{Dudley:2002}.
Combining this with the asymptotic properties of classical group sequential designs, we conclude
\begin{align*}
\rho_K\left(\mathcal{L}\left({S}^\tau|\biguplus_{i=1}^K{Z}_i\right),\mathcal{L}({S})\right)\overset{p}{\rightarrow}0,
\end{align*}
which means that the conditional distribution of ${S}^\tau$ given the observed data approximates the unconditional distribution of ${S}$.
Furthermore, Theorem \ref{theorem-1} states that the limiting multivariate normal distribution of ${S}^\tau$ does not depend on the distribution of the data, particularly, it is achieved for the null hypothesis $H_0$ and the alternative hypothesis $H_1$.
In addition, Theorem \ref{theorem-1} shows that the stage-wise permutation keeps independent increments property for covariance structure which is a favorable property in group sequential designs. See more details about the benefit of independent increments property in \cite{KimTsiatis:2020}.

With the definition of $\varphi_n$ in (\ref{phi}), we give the following definition of its permutated version.
Define the permutated version of $\varphi_n$ as
\begin{align}\label{phi-permutated}
\varphi_n^\tau = \sum_{i=1}^K\varphi_{n,i}^\tau,
\end{align}
where
\begin{equation*}
\varphi_{n,1}^\tau= \left\{
    \begin{aligned}
    1 & \quad \mbox{if} \ S_{\widetilde{m}_1,\widetilde{n}_1}\geq c_1^\tau  \\
    0 & \quad \mbox{otherwise}
    \end{aligned}
    \right.
\end{equation*}
and
\begin{equation*}
\varphi_{n,i}^\tau= \left\{
    \begin{aligned}
    1 & \quad \mbox{if} \ S_{\widetilde{m}_1,\widetilde{n}_1}< c_1^\tau \ \mbox{and} \ \ldots \ \mbox{and} \ S_{\widetilde{m}_{i-1},\widetilde{n}_{i-1}}< c_{i-1}^\tau \ \mbox{and} \ S_{\widetilde{m}_i,\widetilde{n}_i}\geq c_i^\tau  \\
    0 & \quad \mbox{otherwise}
    \end{aligned}
    \right.
\end{equation*}
for $i=2,\ldots,K$.
The critical values $\{c_1^\tau,\ldots,c_K^\tau\}$ in $\varphi_n^\tau$ are defined in (\ref{c-i-tau}).
The main difference between $\varphi_n$ defined in (\ref{phi}) and $\varphi_n^\tau$ defined in (\ref{phi-permutated}) is that the critical value $c_i$ depends on the multivariate normal distribution while the critical value $c_i^\tau$ depends on the permutation distribution of the observed data.

In the next theorem, we will show that the unconditional test $\varphi_n$ and the conditional permutation test $\varphi_n^\tau$ are asymptotically equivalent, which means that both tests have asymptotically the same power.

\begin{theorem}\label{theorem-2}
Suppose that the assumptions of Theorem \ref{theorem-1} are fulfilled.

  (i).
  Under the null hypothesis $H_0: \mu_1=\mu_2$,
  the permutation test $\varphi_n^\tau$ is asymptotically exact at significance level $\alpha$, i.e. ${E}(\varphi_n^\tau)\rightarrow\alpha$, and asymptotically equivalent to $\varphi_n$, i.e.
  \begin{align*}
  {E}(\mid\varphi_n-\varphi_n^\tau\mid)\rightarrow 0, \ \mbox{as} \ n  \rightarrow \infty.
  \end{align*}

  (ii). The permutation test $\varphi_n^\tau$ is consistent, i.e.
  \begin{align*}
  {E}(\varphi_n^\tau) \rightarrow \alpha\mathbb{I}\{\mu_1-\mu_2=0\}+\mathbb{I}\{\mu_1-\mu_2\neq0\}, \ \mbox{as} \ n\rightarrow\infty,
  \end{align*}
  where $\mathbb{I}\{\cdot\}$ is the indicator function.


\end{theorem}

Theorem \ref{theorem-1} and Theorem \ref{theorem-2} show that the studentized permutation test $\varphi_n^\tau$ is asymptotically an appropriate level $\alpha$ test procedure for $H_0: \mu_1=\mu_2$.
In practical applications, it is interesting to explore whether the proposed method provides a more accurate approximation of the distribution of the test statistics than the classical group sequential designs, especially for small sample sizes.
In the next section, we will provide extensive simulation results.

\section{Simulations}\label{sec-4}

We investigate the finite sample properties of the proposed permutation test $\varphi_n^\tau$ in comparison with the classical test $\varphi_n$ and its $t$-approximation $\tilde{\varphi}_n$, with regard to (i) maintaining the Type I error rate under the null hypothesis $H_0: \mu_1=\mu_2$ and (ii) evaluating their power under specific alternative hypothesis $H_1$.

All simulations were conducted using R computing environment, version 4.2.0 \citep{R:2022}, each with 10000 simulations per scenario and 10000 permutations.
Due to the abundance of different group sequential designs, we restrict the analyses to two-, three-, and five-stage designs generating observations from normal distributions as well as a range of non-normal distributions.
The latter is of particular interest from a practical point of view because the assumption of normality may not be satisfied in applications.
Several commonly used non-normal distributions such as the $t$-distribution, the exponential distribution, the Laplace distribution, and the log-normal distribution are considered.

We use $\alpha=2.5\%$ for one-sided hypothesis testing and choose different types of $F_1$ and $F_2$ under the null hypothesis $H_0: \mu_1=\mu_2$.
For example, the scenario $N(0,1)$ vs $N(0,1)$ in the following figures means that $F_1$ represents the standard normal distribution in the treatment arm and $F_2$ represents the standard normal distribution in the control arm.
For non-normal distributions, $F_1$ and $F_2$ could be replaced with the corresponding distributions.
The probability density functions of $F_1$ and $F_2$ selected for the subsequent simulations are depicted in the supplementary material.

\begin{figure}[hpt!]
  \centering
  \includegraphics[width=420pt]{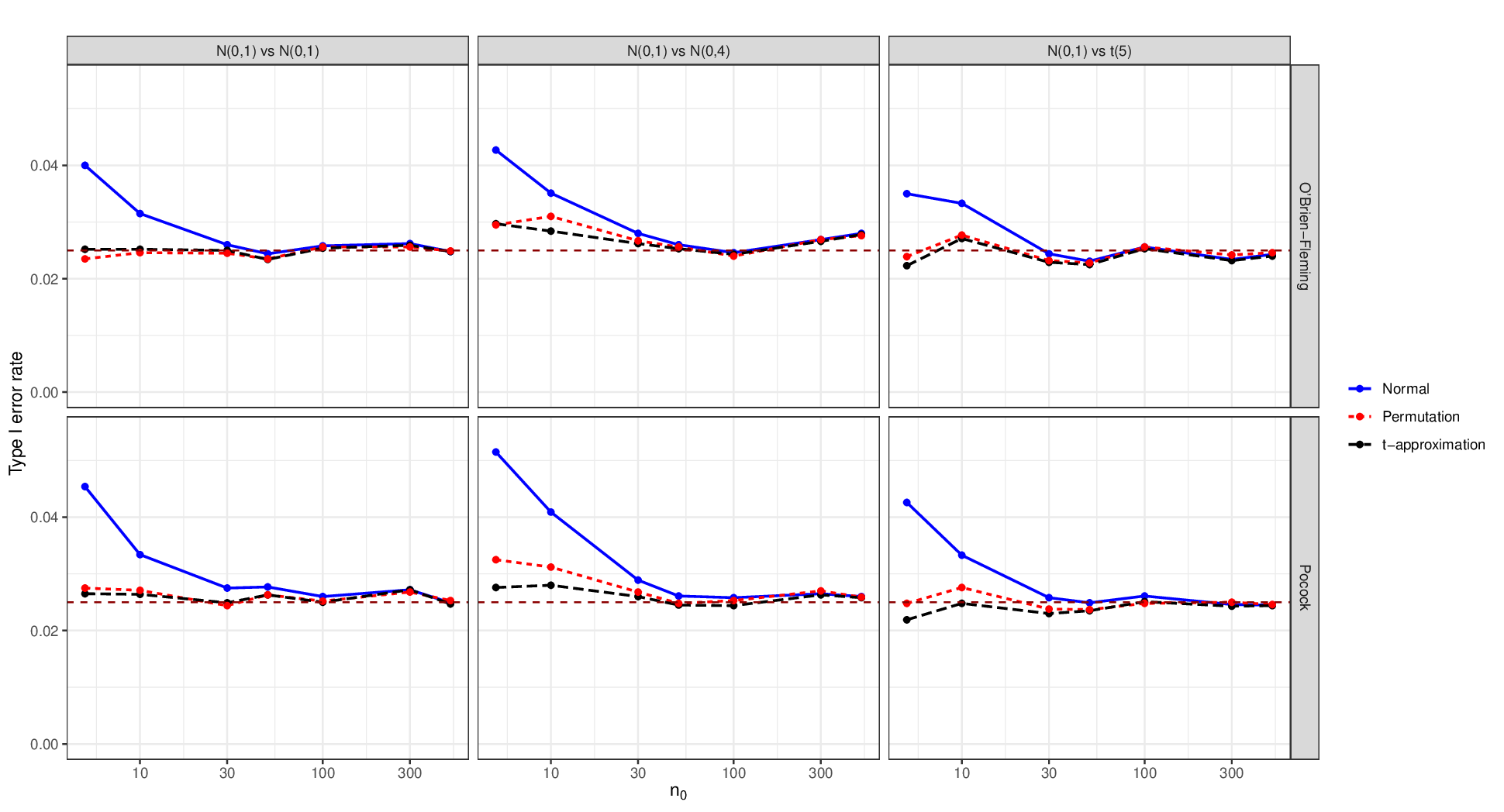}\\
  \caption{Type I error rate ($\alpha=0.025$) simulation results ($y$-axis) of the classical normal method, its $t$-approximation and the proposed permutation method for various scenarios, sample size increments $n_0\in\{5,10,30,50,100,300,500\}$ ($x$-axis), one-sided hypothesis testing, 1:1 allocation ratio and two-stage design}
  \label{two-stages-1}
\end{figure}

\begin{figure}[hpt!]
  \centering
  \includegraphics[width=420pt]{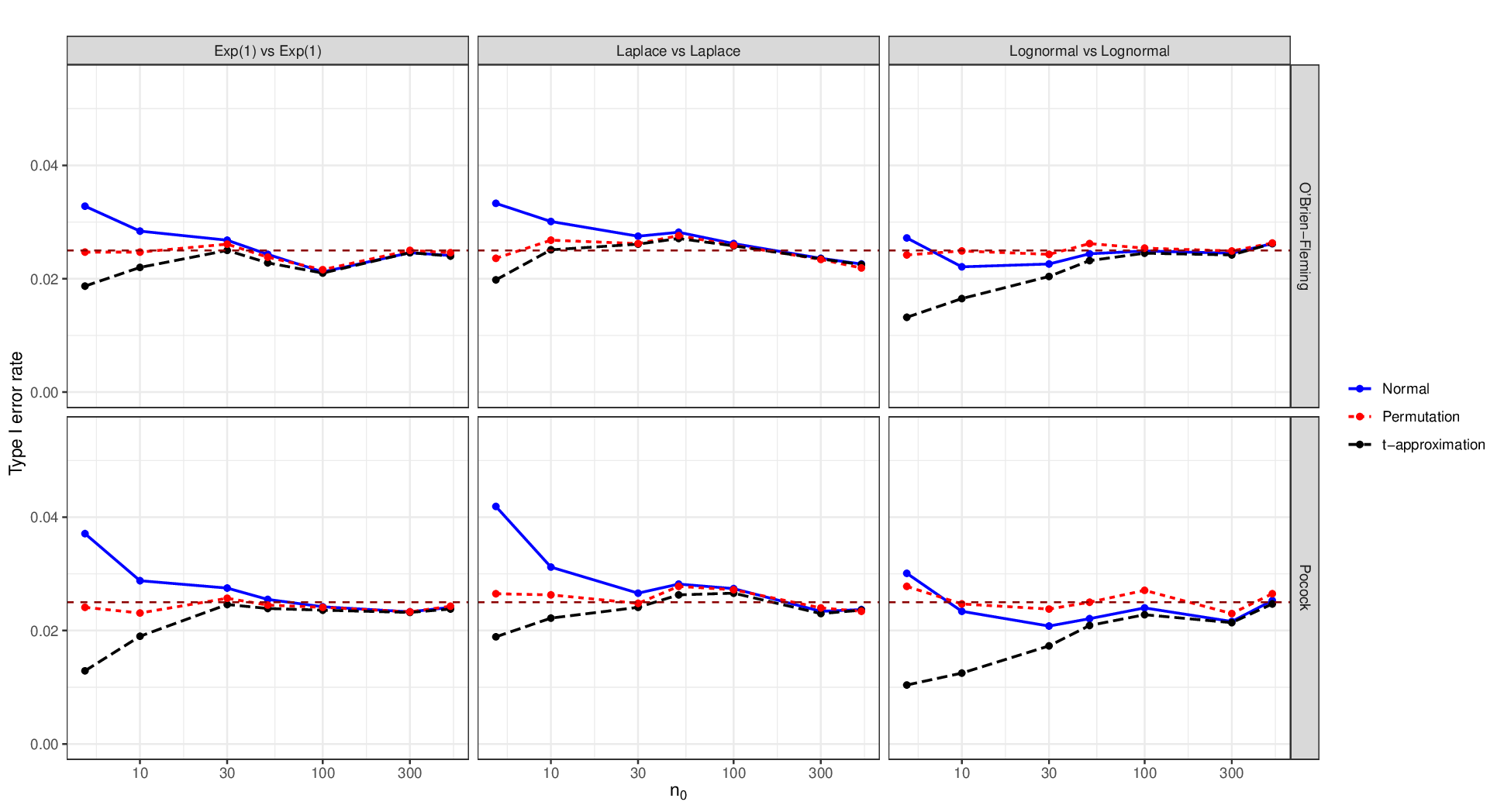}\\
  \caption{Type I error rate ($\alpha=0.025$) simulation results ($y$-axis) of the classical normal method, its $t$-approximation and the proposed permutation method for various scenarios, sample size increments $n_0\in\{5,10,30,50,100,300,500\}$ ($x$-axis), one-sided hypothesis testing, 1:1 allocation ratio and two-stage design}
  \label{two-stages-2}
\end{figure}

In order to simultaneously consider both smaller and larger sample sizes, we set the sample size in the control arm at any given stage, $n_0$, to $5,10,30,50,100,300,500$.
The simulation results for the Type I error rate in balanced designs are displayed in Figure \ref{two-stages-1} and Figure \ref{two-stages-2}.
More detailed simulation results for three-stage and five-stage designs can be found in the supplementary material.
For convenience, these figures are plotted on a logarithmic scale for the $x$-axis.
These results show that the classical group sequential design and its $t$-approximation that was proposed to address small sample problems fail to control the Type I error rate in the case of small sample sizes.
Particularly, the classical test tends to be liberal, while its $t$-approximation tends to be conservative in some scenarios.
It is noted that the $t$-approximation can improve the classical test, but the effect is still limited.
For skewed distributions, the $t$-approximation may yield liberal or conservative results.
With large sample sizes, both the classical test and its $t$-approximation are adequate for the simulated designs.
In contrast, the studentized permutation test controls the Type I error rate very well and greatly improves the existing tests.
For situations with exchangeability between the treatment arm and the control arm, i.e. the observations in the arms follow the same distribution, the studentized permutation test exhibits better control of the Type I error rate than the classical group sequential design and its $t$-approximation.
In addition, the inflation of the Type I error rate by some methods is generally more pronounced for the Pocock type $\alpha$-spending function than for the O’Brien-Fleming type $\alpha$-spending function.

\begin{figure}[hpt!]
  \centering
  \includegraphics[width=420pt]{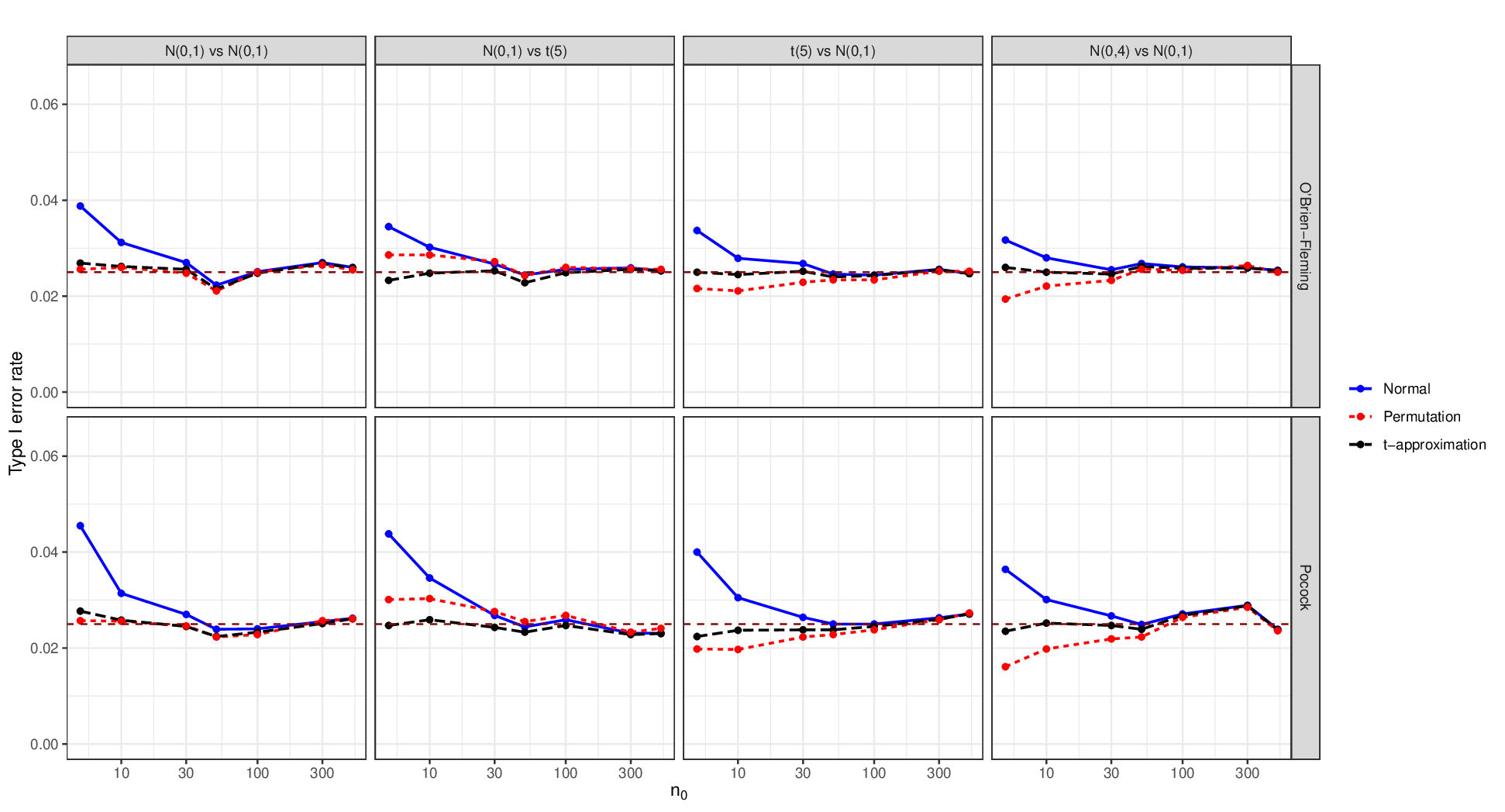}\\
  \caption{Type I error rate ($\alpha=0.025$) simulation results ($y$-axis) of the classical normal method, its $t$-approximation and the proposed permutation method for various scenarios, sample size increments $n_0\in\{5,10,30,50,100,300,500\}$ ($x$-axis), one-sided hypothesis testing, 2:1 allocation ratio and two-stage design}
  \label{unequal-1}
\end{figure}

\begin{figure}[hpt!]
  \centering
  \includegraphics[width=420pt]{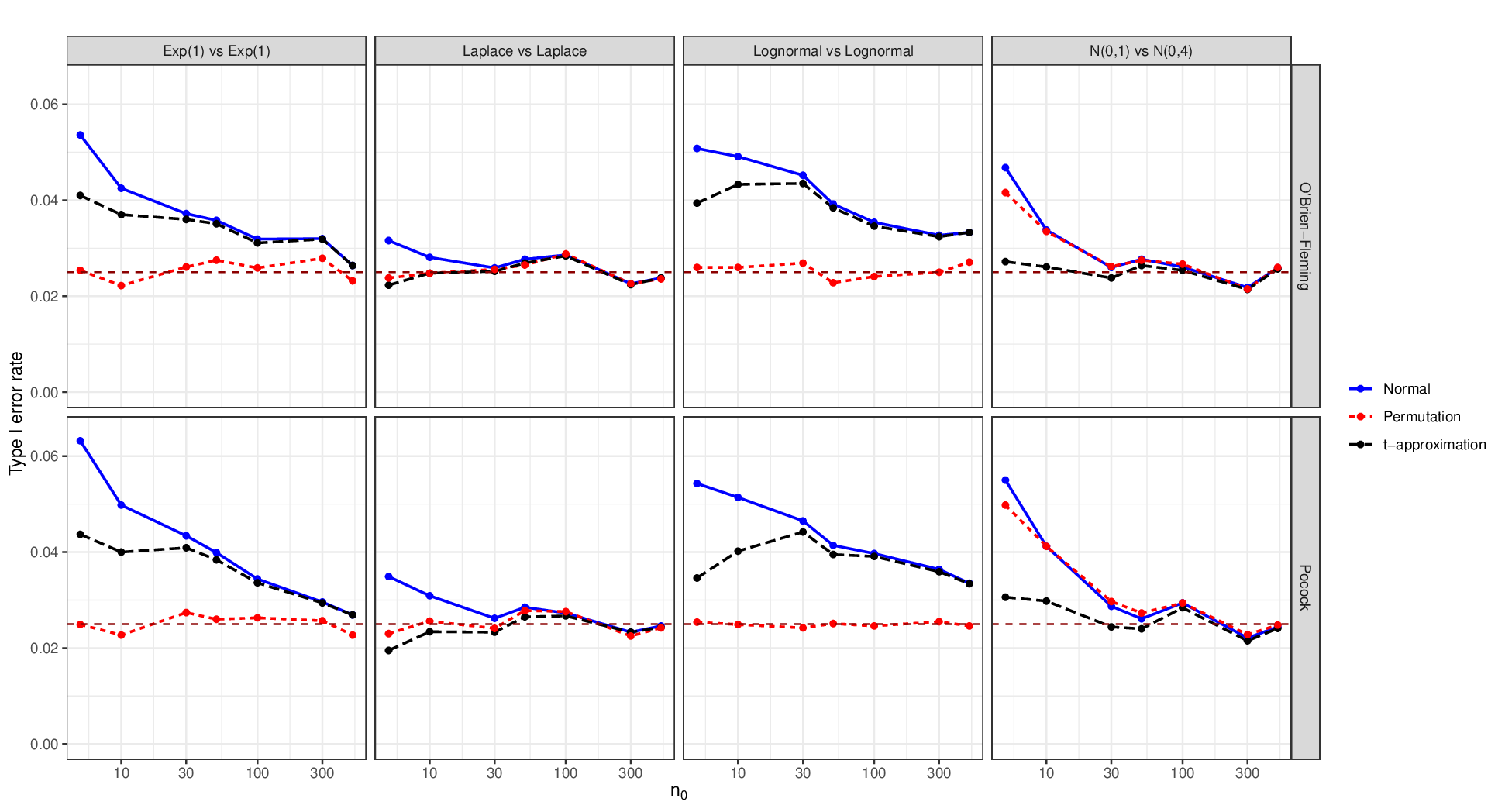}\\
  \caption{Type I error rate ($\alpha=0.025$) simulation results ($y$-axis) of the classical normal method, its $t$-approximation and the proposed permutation method for various scenarios, sample size increments $n_0\in\{5,10,30,50,100,300,500\}$ ($x$-axis), one-sided hypothesis testing, 2:1 allocation ratio and two-stage design}
  \label{unequal-2}
\end{figure}

Next, we examine the influence of the allocation ratio $\gamma$ on the Type I error rate.
In Figure \ref{unequal-1} and Figure \ref{unequal-2}, we present this influence for normal distributions and non-normal distributions.
In unbalanced designs, the $t$-approximation outperforms both the classical test and the proposed test for the normal and $t$-distribution.
This highlights the advantage of $t$-approximation specifically designed for normally or $t$-distributed data, allowing it to overcome heteroscedasticity between the treatment arm and the control arm.
It is important to note that the proposed test may not control the Type I error rate in the case where unequal sample sizes and heteroscedasticity happen between the treatment arm and the control arm. This behaviour is similar to \citet[Table 2]{Janssen:1997} when considering the fixed-sample design.

Finally, we conduct a power comparison among these three methods.
For illustration, we consider the scenario $N(\mu,1)$ vs $N(0,1)$ with varying values of $\mu$.
The simulation results are displayed in Figure \ref{power-1}.
More detailed simulation results for other distributions can be found in the supplementary material.
To illustrate their power, we select $n_0=5, 10, 30$.
Although the classical test does not control the Type I error rate for small samples, we have included the statistic in the power simulation study for illustrative purposes.
These simulations demonstrate that the power curves of the $t$-approximation and the proposed test coincide for small and large sample sizes.
Furthermore, the power curves from the three competing methods coincide asymptotically, suggesting that no differences among these methods emerge with large sample sizes.
Particularly, both the classical test and the proposed permutation test have a comparable power for large sample sizes, which is in line with Theorem \ref{theorem-2}.
Given that the classical test does not adequately control the Type I error rate in small sample settings, we emphasize the efficacy of the proposed permutation test.

In this extensive simulation study, the proposed test controls Type I error rate well for small sample sizes while preserving the same asymptotic properties as classical group sequential designs.
The classical test $\varphi_n$ is valid for large sample sizes; however, it may dramatically fail for small sample sizes.
The $t$-approximation tends to yield conservative results in the case of small sample sizes.
Even under non-exchangeable observations, the proposed test effectively controls the Type I error rate, and this holds true for skewed distributions.

\begin{figure}[hpt!]
  \centering
  \includegraphics[width=420pt]{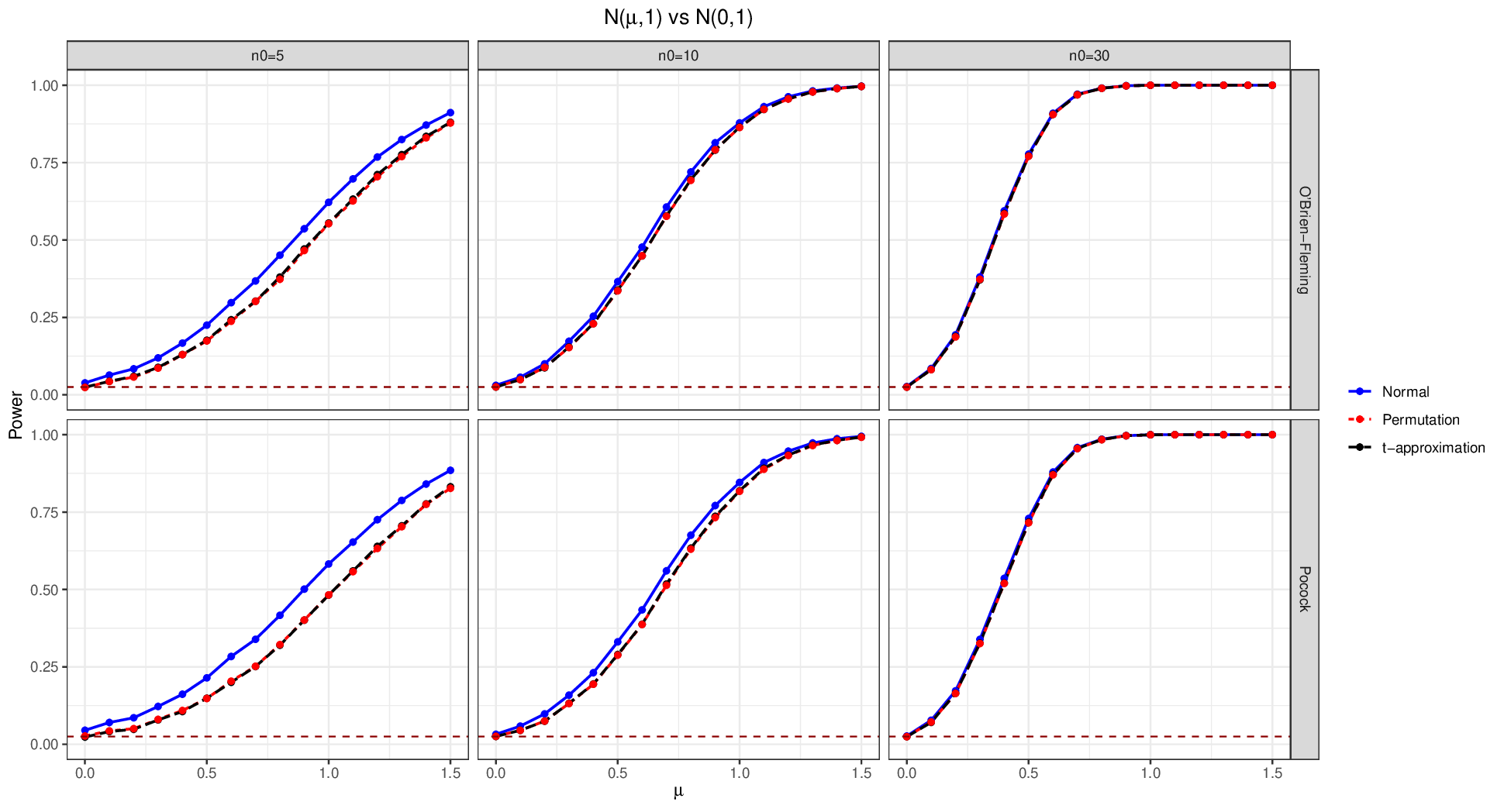}\\
  \caption{Power simulation results ($y$-axis) of the classical normal method, its $t$-approximation and the proposed permutation method for the scenario $N(\mu,1)$ vs $N(0,1)$ under the alternative $\mu\in\{0,0.1,0.2,\ldots,1.5\}$ ($x$-axis), sample size $n_0=5,10,30$, one-sided hypothesis testing, 1:1 allocation ratio and two-stage design}
  \label{power-1}
\end{figure}



\section{Concluding remarks}\label{sec-conclusion}

Based on the ideas introduced by \cite{Janssen:1997} and \cite{ChungRomano:2013}, we have developed a permutation test for group sequential designs including unbalanced designs and non-exchangeable data and explored its finite sample and asymptotic properties.
We have shown that the permutation test shares the same asymptotical properties with the classical group sequential designs. Furthermore, we have shown that both the unconditional test $\varphi_n$ and the conditional permutation test $\varphi_n^\tau$ have asymptotically the same power.
The extensive simulation study illustrates that the proposed test controls Type I error rate well for small sample sizes.

In this paper, we have focused on the comparison between two arms using the test statistic defined in (\ref{eq-statistics}), which is a standard case in medical research.
We hope this work sparks interest in the application of studentized permutation tests in group sequential designs.
Currently, we are exploring the extension to general U-statistics.
The class of U-statistics includes not only commonly-used statistics such as mean and variance but also the Wilcoxon statistic which can be directly applied to the Wilcoxon test without the assumption of the continuity of the underlying distributions.
While the results of fixed-sample design for general U-statistics are discussed in \cite{ChungRomano:2016}, this extension to group sequential design appears to be more difficult because of the consideration of contiguity (cf. \citet[Section 12.3]{LehmannRomano:2005}).
\cite{Nowaketal:2022} studied classical group sequential designs for the Wilcoxon test and indeed observed the inflation of Type I error rate in some settings.
In the future, we plan to apply the studentized permutation test to rank-based statistics in group sequential designs to control the Type I error rate.
We are also interested in extending our method to various data types and hypotheses including both parametric and nonparametric methods.
In addition, the comparison among multiple arms may be a future research direction.

\section*{Acknowledgement}
The authors are grateful for funding by the DFG grant number FR 3070/4-1.

\section*{Supplementary material}
\label{SM}
The supplementary material contains proofs of the theoretical results and additional simulation results.
The R code for the simulation results can be found at \url{https://github.com/xu19971997/a-studentized-permutation-test-in-group-sequential-designs/blob/main/R-code.R}.

\pagebreak
\begin{center}
\textbf{\large Supplementary Material: A studentized permutation test in group sequential designs}
\end{center}
\setcounter{equation}{0}
\setcounter{figure}{0}
\setcounter{table}{0}
\setcounter{page}{1}
\makeatletter
\renewcommand{\theequation}{S\arabic{equation}}
\renewcommand{\thefigure}{S\arabic{figure}}
\renewcommand{\bibnumfmt}[1]{[S#1]}
\renewcommand{\citenumfont}[1]{S#1}

\begin{abstract}
In this supplementary material, we provide proofs of the theoretical results in the main text, along with additional simulation results illustrating the efficacy of the proposed test. In addition, the probability density functions of $F_1$ and $F_2$ selected for the simulations are depicted in Table \ref{Table-1}.
\end{abstract}

\section{The proofs}

Before we provide the proofs of Theorem \ref{theorem-1} and Theorem \ref{theorem-2}, we need the following two important limits (\ref{limit-1}) and (\ref{limit-2}).
Note that
\begin{align}\label{eq-var}
\mbox{var}\left(\widehat{\mu}_{1,k}-\widehat{\mu}_{2,k}\right)=\mbox{var}\left(\frac{1}{\widetilde{m}_k}\sum_{i=1}^{\widetilde{m}_k}X_i-\frac{1}{\widetilde{n}_k}\sum_{j=1}^{\widetilde{n}_k}Y_j\right)=\frac{\sigma_1^2}{\widetilde{m}_k}+\frac{\sigma_2^2}{\widetilde{n}_k}=\left(\frac{\sigma_1^2}{\gamma}+\sigma_2^2\right)\frac{1}{\widetilde{n}_k}.
\end{align}
Substituting (\ref{eq-var}) into (\ref{assumption-GSD}), we have
\begin{align*}
\lim_{n\rightarrow\infty}\frac{\widetilde{n}_k}{n}=\left(\frac{\sigma_1^2}{\gamma}+\sigma_2^2\right)\mathcal{I}_k^*.
\end{align*}

Therefore, we obtain the two important limits
\begin{align}\label{limit-1}
\frac{\sum_{k=1}^{\min\{i,j\}} n_k}{\sum_{k=1}^{\max\{i,j\}} n_k}=\frac{\widetilde{n}_{\min\{i,j\}}}{\widetilde{n}_{\max\{i,j\}}}=\frac{\widetilde{n}_{\min\{i,j\}}/n}{\widetilde{n}_{\max\{i,j\}}/n}\rightarrow \frac{\mathcal{I}_{\min\{i,j\}}^*}{\mathcal{I}_{\max\{i,j\}}^*} \quad \mbox{as} \ n\rightarrow\infty,
\end{align}
and
\begin{align}\label{limit-2}
\frac{n_j}{\sum_{i=1}^k n_i} = \frac{\widetilde{n}_j-\widetilde{n}_{j-1}}{\widetilde{n}_k} = \frac{\widetilde{n}_j/n-\widetilde{n}_{j-1}/n}{\widetilde{n}_k/n} \rightarrow \frac{\mathcal{I}_j^*-\mathcal{I}_{j-1}^*}{\mathcal{I}_k^*} \quad \mbox{as} \ n\rightarrow\infty.
\end{align}
For convenience, set $\mathcal{I}_0^*=0$ in (\ref{limit-2}).

To prove Theorem 1, we first establish the useful lemma.

\begin{lemma}\label{lemma-2}
For all $k=1,\ldots,K$, given the observed data $\biguplus_{i=1}^k {Z}_i$, we have
\begin{align*}
\widehat{\sigma}_{1,k,\tau}^2 \overset{p}{\rightarrow} \sigma_0^2 \ \mbox{and} \ \widehat{\sigma}_{2,k,\tau}^2 \overset{p}{\rightarrow} \sigma_0^2,
\end{align*}
where $\sigma_0^2=\mbox{var}_{X'\sim H}(X')$ with the mixture distribution $H=\frac{\gamma}{\gamma+1}F_1+\frac{1}{\gamma+1}F_2$.
\end{lemma}

\noindent\textit{Proof of Lemma \ref{lemma-2}.}

Due to the similar expression of $\widehat{\sigma}_{1,k,\tau}^2$ and $\widehat{\sigma}_{2,k,\tau}^2$, it is sufficient to prove $\widehat{\sigma}_{1,k,\tau}^2 \overset{p}{\rightarrow} \sigma_0^2$.

First, we consider, for a specific $u\in\{1,\ldots,k\}$,
\begin{align*}
&\sum_{i=1}^{m_u}\left(Z_{u,\tau_u(i)}-\frac{1}{\widetilde{m}_k}\sum_{s=1}^{k}\sum_{j=1}^{m_s}Z_{s,\tau_s(j)}\right)^2 \\
&= \sum_{i=1}^{m_u}\left(Z_{u,\tau_u(i)}-\frac{1}{m_u}\sum_{j=1}^{m_u}Z_{u,\tau_u(j)}+\underbrace{\frac{1}{m_u}\sum_{j=1}^{m_u}Z_{u,\tau_u(j)}-\frac{1}{\widetilde{m}_k}\sum_{j=1}^{m_u}Z_{u,\tau_u(j)}}_{A_u}-\underbrace{\frac{1}{\widetilde{m}_k}\sum_{s\neq u}\sum_{j=1}^{m_s}Z_{s,\tau_s(j)}}_{B_u}\right)^2\\
&= \sum_{i=1}^{m_u}\left(Z_{u,\tau_u(i)}-\frac{1}{m_u}\sum_{j=1}^{m_u}Z_{u,\tau_u(j)}+A_u-B_u\right)^2\\
&= \sum_{i=1}^{m_u}\left\{ \left(Z_{u,\tau_u(i)}-\frac{1}{m_u}\sum_{j=1}^{m_u}Z_{u,\tau_u(j)}\right)^2 + 2\left(Z_{u,\tau_u(i)}-\frac{1}{m_u}\sum_{j=1}^{m_u}Z_{u,\tau_u(j)}\right)(A_u-B_u) + (A_u-B_u)^2\right\}\\
&= \sum_{i=1}^{m_u}\left(Z_{u,\tau_u(i)}-\frac{1}{m_u}\sum_{j=1}^{m_u}Z_{u,\tau_u(j)}\right)^2 + \sum_{i=1}^{m_u}2\left(Z_{u,\tau_u(i)}-\frac{1}{m_u}\sum_{j=1}^{m_u}Z_{u,\tau_u(j)}\right)(A_u-B_u) + \sum_{i=1}^{m_u}(A_u-B_u)^2\\
&= \sum_{i=1}^{m_u}\left(Z_{u,\tau_u(i)}-\frac{1}{m_u}\sum_{j=1}^{m_u}Z_{u,\tau_u(j)}\right)^2 + m_u(A_u-B_u)^2.
\end{align*}

Thus, we rewrite
\begin{align*}
\widehat{\sigma}_{1,k,\tau}^2&= \frac{1}{\widetilde{m}_k-1}\sum_{u=1}^{k}\sum_{i=1}^{m_u}\left(Z_{u,\tau_u(i)}-\frac{1}{\widetilde{m}_k}\sum_{s=1}^{k}\sum_{j=1}^{m_s}Z_{s,\tau_s(j)}\right)^2\\
&= \frac{1}{\widetilde{m}_k-1}\sum_{u=1}^{k}\left\{\sum_{i=1}^{m_u}\left(Z_{u,\tau_u(i)}-\frac{1}{m_u}\sum_{j=1}^{m_u}Z_{u,\tau_u(j)}\right)^2 + m_u(A_u-B_u)^2\right\}\\
&= \frac{\widetilde{m}_k}{\widetilde{m}_k-1}\frac{1}{\sum_{i=1}^{k}m_i}\sum_{u=1}^{k}\left\{m_u\frac{1}{m_u}\sum_{i=1}^{m_u}\left(Z_{u,\tau_u(i)}-\frac{1}{m_u}\sum_{j=1}^{m_u}Z_{u,\tau_u(j)}\right)^2 + m_u(A_u-B_u)^2\right\}\\
&= \frac{\widetilde{m}_k}{\widetilde{m}_k-1}\left\{\sum_{u=1}^{k}\frac{m_u}{\widetilde{m}_k}\frac{1}{m_u}\sum_{i=1}^{m_u}\left(Z_{u,\tau_u(i)}-\frac{1}{m_u}\sum_{j=1}^{m_u}Z_{u,\tau_u(j)}\right)^2+\sum_{u=1}^{k}\frac{m_u}{\widetilde{m}_k}(A_u-B_u)^2\right\}.
\end{align*}

Applying Theorem 3.7.1. in \cite{VaartWellner:1996}, we have
\begin{align}\label{eq-V-n}
\frac{1}{m_u}\sum_{i=1}^{m_u}\left(Z_{u,\tau_u(i)}-\frac{1}{m_u}\sum_{j=1}^{m_u}Z_{u,\tau_u(j)}\right)^2 \overset{p}{\rightarrow}\sigma_0^2 , \ u=1,\ldots,k,
\end{align}
and
\begin{align}\label{eq-U-n}
\frac{1}{m_u}\sum_{j=1}^{m_u}Z_{u,\tau_u(j)}-\frac{1}{m_s}\sum_{j=1}^{m_s}Z_{s,\tau_s(j)} \overset{p}{\rightarrow}0 , \ s\neq u.
\end{align}

Notice that we can rewrite $A_u$ and $B_u$ as
\begin{align*}
A_u &= \sum_{s\neq u}\frac{m_s}{\widetilde{m}_k}\frac{1}{m_u}\sum_{j=1}^{m_u}Z_{u,\tau_u(j)},\\
B_u &= \sum_{s\neq u}\frac{m_s}{\widetilde{m}_k}\frac{1}{m_s}\sum_{j=1}^{m_s}Z_{s,\tau_s(j)}.
\end{align*}

Thus, combined with (\ref{limit-2}) and (\ref{eq-U-n}), we have
\begin{align}\label{Au-Bu}
A_u-B_u &= \sum_{s\neq u}\frac{m_s}{\widetilde{m}_k}\left( \frac{1}{m_u}\sum_{j=1}^{m_u}Z_{u,\tau_u(j)}-\frac{1}{m_s}\sum_{j=1}^{m_s}Z_{s,\tau_s(j)} \right) \overset{p}{\rightarrow} 0 .
\end{align}

Combined (\ref{limit-2}), (\ref{eq-V-n}), and (\ref{Au-Bu}) with the expression of $\widehat{\sigma}_{1,k,\tau}^2$, we have $\widehat{\sigma}_{1,k,\tau}^2 \overset{p}{\rightarrow} \sigma_0^2$, which completes this proof.

\begin{flushright}
$\square$
\end{flushright}

\noindent\textit{Proof of Theorem 1.}

Define
\begin{align*}
\widehat{\sigma}^2_{0,k,\tau} &= \frac{\widetilde{n}_k}{\widetilde{m}_k+\widetilde{n}_k}\widehat{\sigma}^2_{1,k,\tau}+\frac{\widetilde{m}_k}{\widetilde{m}_k+\widetilde{n}_k}\widehat{\sigma}^2_{2,k,\tau}, \ k=1,\ldots,K.
\end{align*}

Since ${m_k}/{n_k}=\gamma$ from (\ref{gamma}), it follows that
\begin{align*}
\widehat{\sigma}^2_{0,k,\tau} &= \frac{1}{\gamma+1}\widehat{\sigma}^2_{1,k,\tau}+\frac{\gamma}{\gamma+1}\widehat{\sigma}^2_{2,k,\tau}.
\end{align*}

We note that ${S}^\tau$ can be decomposed as
\begin{align*}
{S}^\tau = {A}^\tau {E}^\tau
\end{align*}
where
\begin{align*}
 {A}^\tau = diag\left(\widehat{\sigma}^{-1}_{0,1,\tau},\ldots,\widehat{\sigma}^{-1}_{0,K,\tau}\right)
\end{align*}
and
\begin{align*}
{E}^\tau = \left(\left(\frac{\widetilde{m}_1\widetilde{n}_1}{\widetilde{m}_1+\widetilde{n}_1}\right)^{1/2}\left(\widehat{\mu}_{1,1}^\tau-\widehat{\mu}_{2,1}^\tau\right),\ldots,\left(\frac{\widetilde{m}_K\widetilde{n}_K}{\widetilde{m}_K+\widetilde{n}_K}\right)^{1/2}\left(\widehat{\mu}_{1,K}^\tau-\widehat{\mu}_{2,K}^\tau\right)\right)^{\T}.
\end{align*}

First, we consider the asymptotic property of ${A}^\tau$.
Due to Lemma \ref{lemma-2}, we know that $\widehat{\sigma}^2_{0,k,\tau}$ converges in probability to $\sigma_0^2$ for $k=1,\ldots,K$. Thus, we have convergence in probability
\begin{align*}
{A}^\tau \rightarrow diag(\sigma_0^{-1},\ldots,\sigma_0^{-1})
\end{align*}
as $n\rightarrow\infty$.

Second, we consider the asymptotic property of ${E}^\tau$.
Similar to \cite{Paulyetal:2015}, we note that the classical Cram\'{e}r-Wold device cannot be applied directly in the context of conditional weak convergence results due to the occurrence of uncountably many exceptional sets.
Therefore, we apply a modified Cram\'{e}r-Wold device (see the proof of Theorem 4.1 in \cite{Pauly:2011}).

For convenience, define
\begin{align*}
D_i = \left({m_i+n_i}\right)^{1/2}\left(\frac{1}{m_i}\sum_{j=1}^{m_i}Z_{i,\tau_i(j)}-\frac{1}{n_i}\sum_{j=m_i+1}^{m_i+n_i}Z_{i,\tau_i(j)}\right),\ i=1,\ldots,K.
\end{align*}
Let $C$ be a dense and countable subset of $\mathbb{R}^K$.
Then for every fixed $
{\lambda}=(\lambda_1,\ldots,\lambda_K)^\T\in C$ we have
\begin{align*}
{\lambda}^\T{E}^\tau &= \sum_{k=1}^K\left\{\left(\frac{\widetilde{m}_k\widetilde{n}_k}{\widetilde{m}_k+\widetilde{n}_k}\right)^{1/2}\sum_{i=1}^k\frac{n_i}{\widetilde{n}_k}\frac{1}{(m_i+n_i)^{1/2}}D_i\right\}\lambda_k\\
&= \sum_{k=1}^K\sum_{i=1}^k\left(\frac{\widetilde{m}_k\widetilde{n}_k}{\widetilde{m}_k+\widetilde{n}_k}\right)^{1/2}\frac{n_i}{\widetilde{n}_k}\frac{1}{(m_i+n_i)^{1/2}}D_i\lambda_k\\
&= \sum_{i=1}^K\sum_{k=i}^K\left(\frac{\widetilde{m}_k\widetilde{n}_k}{\widetilde{m}_k+\widetilde{n}_k}\right)^{1/2}\frac{n_i}{\widetilde{n}_k}\frac{1}{(m_i+n_i)^{1/2}}D_i\lambda_k \ (\mbox{since} \ \sum_{i=1}^n\sum_{j=1}^i a_ib_j=\sum_{j=1}^n\sum_{i=j}^n a_ib_j)\\
&= \sum_{i=1}^K\left\{\sum_{k=i}^K\left(\frac{\widetilde{m}_k\widetilde{n}_k}{\widetilde{m}_k+\widetilde{n}_k}\right)^{1/2}\frac{\lambda_k}{\widetilde{n}_k}\right\}\frac{n_i}{(m_i+n_i)^{1/2}}D_i.
\end{align*}

From Theorem 2.1. in \cite{ChungRomano:2013}, we know
\begin{align*}
D_i \overset{d}{\rightarrow} N\left(0,\frac{(\gamma+1)^2}{\gamma}\sigma_0^2\right), \ i=1,\ldots,K,
\end{align*}
in probability.
Note that $D_i$ are independent.
Altogether with (\ref{limit-2}), this implies by Slutsky's theorem
\begin{align}\label{conv}
{\lambda}^\T{E}^\tau \overset{d}{\rightarrow} N(0,\sigma_{\lambda}^2)
\end{align}
in probability, where
\begin{align*}
\sigma_{\lambda}^2 &= \frac{\gamma}{(\gamma+1)^2} \sum_{i=1}^K \left\{ \sum_{k=i}^K\lambda_k\left(\frac{\mathcal{I}_i^*-\mathcal{I}_{i-1}^*}{\mathcal{I}_k^*}\right)^{1/2} \right\}^2 \frac{(\gamma+1)^2}{\gamma}\sigma_0^2 \\
&= \sigma_0^2 \sum_{i=1}^K \left\{ \sum_{k=i}^K\lambda_k\left(\frac{\mathcal{I}_i^*-\mathcal{I}_{i-1}^*}{\mathcal{I}_k^*}\right)^{1/2} \right\}^2 \\
&= \sigma_0^2 \sum_{i=1}^K (\mathcal{I}_i^*-\mathcal{I}_{i-1}^*) \left\{ \sum_{k=i}^K\frac{\lambda_k}{(\mathcal{I}_k^*)^{1/2}} \right\}^2 \\
&= \sigma_0^2 \sum_{i=1}^K (\mathcal{I}_i^*-\mathcal{I}_{i-1}^*) \left\{ \sum_{k=i}^K\frac{\lambda_k^2}{\mathcal{I}_k^*} + \sum_{i\leq i'< j'\leq K} \frac{\lambda_{i'}}{(\mathcal{I}_{i'}^*)^{1/2}}\frac{\lambda_{j'}}{(\mathcal{I}_{j'}^*)^{1/2}}  \right\} \\
&= \sigma_0^2 \left\{ \sum_{i=1}^K\sum_{k=i}^K (\mathcal{I}_i^*-\mathcal{I}_{i-1}^*)\frac{\lambda_k^2}{\mathcal{I}_k^*} + \sum_{i=1}^K\sum_{i\leq i'< j'\leq K}(\mathcal{I}_i^*-\mathcal{I}_{i-1}^*)\frac{\lambda_{i'}}{(\mathcal{I}_{i'}^*)^{1/2}}\frac{\lambda_{j'}}{(\mathcal{I}_{j'}^*)^{1/2}}  \right\} \\
&= \sigma_0^2 \left\{ \sum_{k=1}^K  \sum_{i=1}^k (\mathcal{I}_i^*-\mathcal{I}_{i-1}^*)\frac{\lambda_k^2}{\mathcal{I}_k^*} + \sum_{1\leq i'< j'\leq K} \sum_{i=1}^{i'}(\mathcal{I}_i^*-\mathcal{I}_{i-1}^*)\frac{\lambda_{i'}}{(\mathcal{I}_{i'}^*)^{1/2}}\frac{\lambda_{j'}}{(\mathcal{I}_{j'}^*)^{1/2}}  \right\} \\
&(\mbox{since} \ \sum_{i=1}^n\sum_{j=1}^i a_ib_j=\sum_{j=1}^n\sum_{i=j}^n a_ib_j, \ \mbox{and} \ \sum_{i=1}^n\sum_{i\leq i'<j'\leq n}a_i b_{i'} c_{j'}=\sum_{1\leq i' < j' \leq n}\sum_{i=1}^{i'}a_i b_{i'} c_{j'}) \\
&= \sigma_0^2 \left\{ \sum_{k=1}^K \lambda_k^2 + \sum_{1\leq i'< j'\leq K} \frac{(\mathcal{I}_{i'}^*)^{1/2}}{(\mathcal{I}_{j'}^*)^{1/2}}\lambda_{i'}\lambda_{j'}  \right\}.
\end{align*}

Note that ${\lambda}^{\T}{U} \sim N(0,\sigma_{\lambda}^2)$
if ${U}=(U_1,\ldots,U_K)^{\T}$ follows a multivariate normal distribution $N({0}_K,{\Sigma}_0)$ where ${0}_K$ denotes the $K$-dimensional vector of zeros and the $(i,j)$th element of covariance matrix $\Sigma_0$ is $\sigma_0^2({\mathcal{I}_{\min\{i,j\}}^*}/{\mathcal{I}_{\max\{i,j\}}^*})^{1/2}$.


Next, we continue to follow the ideas in \cite{Paulyetal:2015}. We note that the convergence (\ref{conv}) holds for every fixed ${\lambda}\in C$.
Applying the subsequential principle for convergence in probability (see e.g. Theorem 9.2.1 in \cite{Dudley:2002}) we can find a common subsequence such that (\ref{conv}) holds almost surely for all ${\lambda}\in C$ along this subsequence.
Now continuity of the characteristic function of the limit and tightness of the conditional distribution of ${E}^\tau$ given the observed data $\biguplus_{i=1}^K{Z}_i$ show that (\ref{conv}) holds almost surely for all ${\lambda}\in\mathbb{R}^K$ along this subsequence.

Thus, an application of the classical Cram\'{e}r-Wold device together with another application of the subsequence principle imply that, the conditional permutation distribution of ${E}^\tau$ given the observed data $\biguplus_{i=1}^K{Z}_i$ weakly converges to a multivariate normal distribution $N({0}_K,{\Sigma}_0)$ in probability.

Combining the asymptotic property of ${A}^\tau$ with ${E}^\tau$, and applying the continuous mapping theorem and Slutsky's theorem, the proof is completed.

\begin{flushright}
$\square$
\end{flushright}

\noindent\textit{Proof of Theorem 2.}

For the first part, we use similar techniques illustrated in \citet[Page 58]{WittingNoelle:1970} (also see Lemma 1 in \cite{JanssenPauls:2003}). Let
\begin{align*}
S_{\widetilde{m}_k,\widetilde{n}_k}^* = S_{\widetilde{m}_k,\widetilde{n}_k} - (c_k^\tau-c_k), \ k=1,\ldots,K.
\end{align*}

Note that (\ref{s-tau-convergence}) implies that the data-dependent critical values $c_i^\tau$ defined in (\ref{c-i-tau}) still converge in probability to $c_i$ defined in (\ref{c-i}) for $i=1,\ldots,K$.
For each $k=1,\ldots,K$, since
\begin{align*}
 \mbox{pr}(S_{\widetilde{m}_k,\widetilde{n}_k}\leq c_k \leq S_{\widetilde{m}_k,\widetilde{n}_k}^*) & = \mbox{pr}(S_{\widetilde{m}_k,\widetilde{n}_k}\leq c_k \leq S_{\widetilde{m}_k,\widetilde{n}_k}^*, \mid S_{\widetilde{m}_k,\widetilde{n}_k}-S_{\widetilde{m}_k,\widetilde{n}_k}^* \mid > \delta(\varepsilon)) \\
& \quad + \mbox{pr}(S_{\widetilde{m}_k,\widetilde{n}_k}\leq c_k \leq S_{\widetilde{m}_k,\widetilde{n}_k}^*, \mid S_{\widetilde{m}_k,\widetilde{n}_k}-S_{\widetilde{m}_k,\widetilde{n}_k}^* \mid \leq \delta(\varepsilon)) \\
& \leq \mbox{pr}(\mid S_{\widetilde{m}_k,\widetilde{n}_k}-S_{\widetilde{m}_k,\widetilde{n}_k}^* \mid > \delta(\varepsilon)) \\
& \quad + \mbox{pr}(S_{\widetilde{m}_k,\widetilde{n}_k}\leq c_k \leq S_{\widetilde{m}_k,\widetilde{n}_k}^*,   S_{\widetilde{m}_k,\widetilde{n}_k}^*-\delta(\varepsilon) \leq S_{\widetilde{m}_k,\widetilde{n}_k} \leq S_{\widetilde{m}_k,\widetilde{n}_k}^* + \delta(\varepsilon)) \\
& \leq \mbox{pr}(\mid c_k^\tau-c_k \mid > \delta(\varepsilon)) + \mbox{pr}(c_k-\delta(\varepsilon) \leq S_{\widetilde{m}_k,\widetilde{n}_k} \leq c_k) \\
& \leq 2\varepsilon \ \mbox{for any} \ n \geq n(\varepsilon),
\end{align*}
it follows that
\begin{align*}
{E}(\mid\varphi_{n,k}-\varphi_{n,k}^\tau\mid) & \leq \mbox{pr}(\mid\varphi_{n,k}-\varphi_{n,k}^\tau\mid > 0) \\
& \leq \mbox{pr}(S_{\widetilde{m}_k,\widetilde{n}_k}\leq c_k \leq S_{\widetilde{m}_k,\widetilde{n}_k}^*)+ \mbox{pr}(S_{\widetilde{m}_k,\widetilde{n}_k}^*\leq c_k \leq S_{\widetilde{m}_k,\widetilde{n}_k})\rightarrow 0, \ \mbox{as} \ n \ \rightarrow \infty.
\end{align*}

Thus, as $n\rightarrow\infty$ we have
\begin{align*}
{E}(\mid\varphi_n-\varphi_n^\tau\mid) = {E}(\mid\sum_{i=1}^K(\varphi_{n,i}-\varphi_{n,i}^\tau)\mid) \leq \sum_{i=1}^K{E}(\mid\varphi_{n,i}-\varphi_{n,i}^\tau\mid) \rightarrow 0.
\end{align*}


With the result ${E}(\varphi_{n})\rightarrow\alpha$ under the null hypothesis $H_0:\mu_1=\mu_2$, the above results imply ${E}(\varphi_{n}^\tau)\rightarrow\alpha$ under the same condition, which completes the first part of this theorem.

Suppose now that $\mu_1\neq\mu_2$.
We consider the numerator of $S_{\widetilde{m}_k,\widetilde{n}_k}$ for $k=1,\ldots,K$. For convenience, let
\begin{align*}
T_{\widetilde{m}_k,\widetilde{n}_k} &= \left(\frac{\widetilde{m}_k\widetilde{n}_k}{\widetilde{m}_k+\widetilde{n}_k}\right)^{1/2}\left(\widehat{\mu}_{1,k}-\widehat{\mu}_{2,k}\right), \ k=1,\ldots,K.
\end{align*}
Since the divergence of the test statistic $T_{\widetilde{m}_k,\widetilde{n}_k}$ under the alternative, it follows that there exists one index $j$ such that
\begin{align*}
T_{\widetilde{m}_j,\widetilde{n}_j} \rightarrow \mbox{sign}(\mu_1-\mu_2)\cdot\infty
\end{align*}
as $n\rightarrow\infty$.
For the denominator of $S_{\widetilde{m}_k,\widetilde{n}_k}$, we know \begin{align*}\frac{\widetilde{n}_k}{\widetilde{m}_k+\widetilde{n}_k}\widehat{\sigma}^2_{1,k}+\frac{\widetilde{m}_k}{\widetilde{m}_k+\widetilde{n}_k}\widehat{\sigma}^2_{2,k}\overset{p}{\rightarrow}\frac{1}{\gamma+1}\sigma^2_{1}+\frac{\gamma}{\gamma+1}\sigma^2_{2}
\end{align*}
for $k=1,\ldots,K$.
This completes the second part of this theorem.

\begin{flushright}
$\square$
\end{flushright}

\section{Additional simulation results}

In this section, we display additional simulation results for Type I error rate and power.
As probability density functions of $F_1$ and $F_2$, we take the functions shown in Table \ref{Table-1}.

\begin{table}
\centering
\caption{Probability density functions}
\label{Table-1}
\begin{tabular}{c|c}
$N(\mu,\sigma^2)$ & $g(x) = {\frac{1}{\sigma{({2\pi})^{1/2}}}}\exp\{-({x-\mu})^{2}/(2{\sigma}^{2})\}$, $x\in \mathbb{R}$ \\
$t(\nu)$ &  $g(x) = \frac{\Gamma\{{({\nu+1})/{2}}\}}{{({\pi\nu})^{1/2}}\Gamma({{\nu}/{2}})}(1+{{x^{2}}/{\nu}})^{-{({\nu+1})/{2}}}$, $x\in \mathbb{R}$ \\
Exp(1) & $g(x) = \exp(-x)$, $x\in[0,\infty)$ \\
Laplace & $g(x) = {\frac{1}{2}}\exp \left(-{|x|}\right)$, $x\in \mathbb{R}$  \\
Log-normal & $g(x) = {\frac{1}{x{({2\pi})^{1/2}}}}\exp\{-{ {(\log x)^{2}}/2}\}$, $x\in(0,\infty)$ \\
\end{tabular}
\end{table}


\subsection{Type I error rate}

\begin{figure}[hpt!]
  \centering
  \includegraphics[width=420pt]{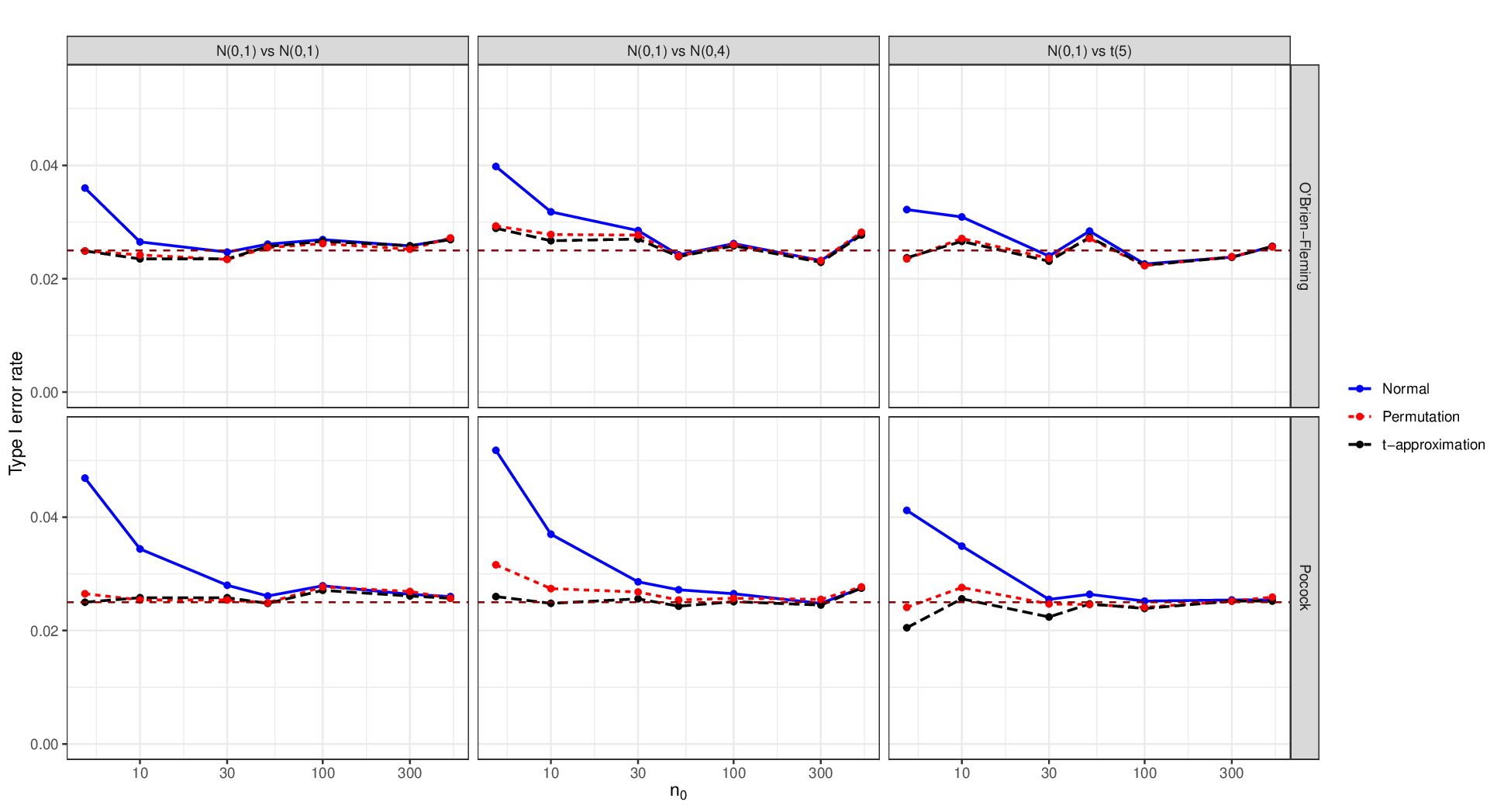}\\
  \caption{Type I error rate ($\alpha=0.025$) simulation results ($y$-axis) of the classical normal method, its $t$-approximation and the proposed permutation method for various scenarios, sample size increments $n_0\in\{5,10,30,50,100,300,500\}$ ($x$-axis), one-sided hypothesis testing, 1:1 allocation ratio and three-stage design}
  \label{three-stages-1}
\end{figure}

\begin{figure}[hpt!]
  \centering
  \includegraphics[width=420pt]{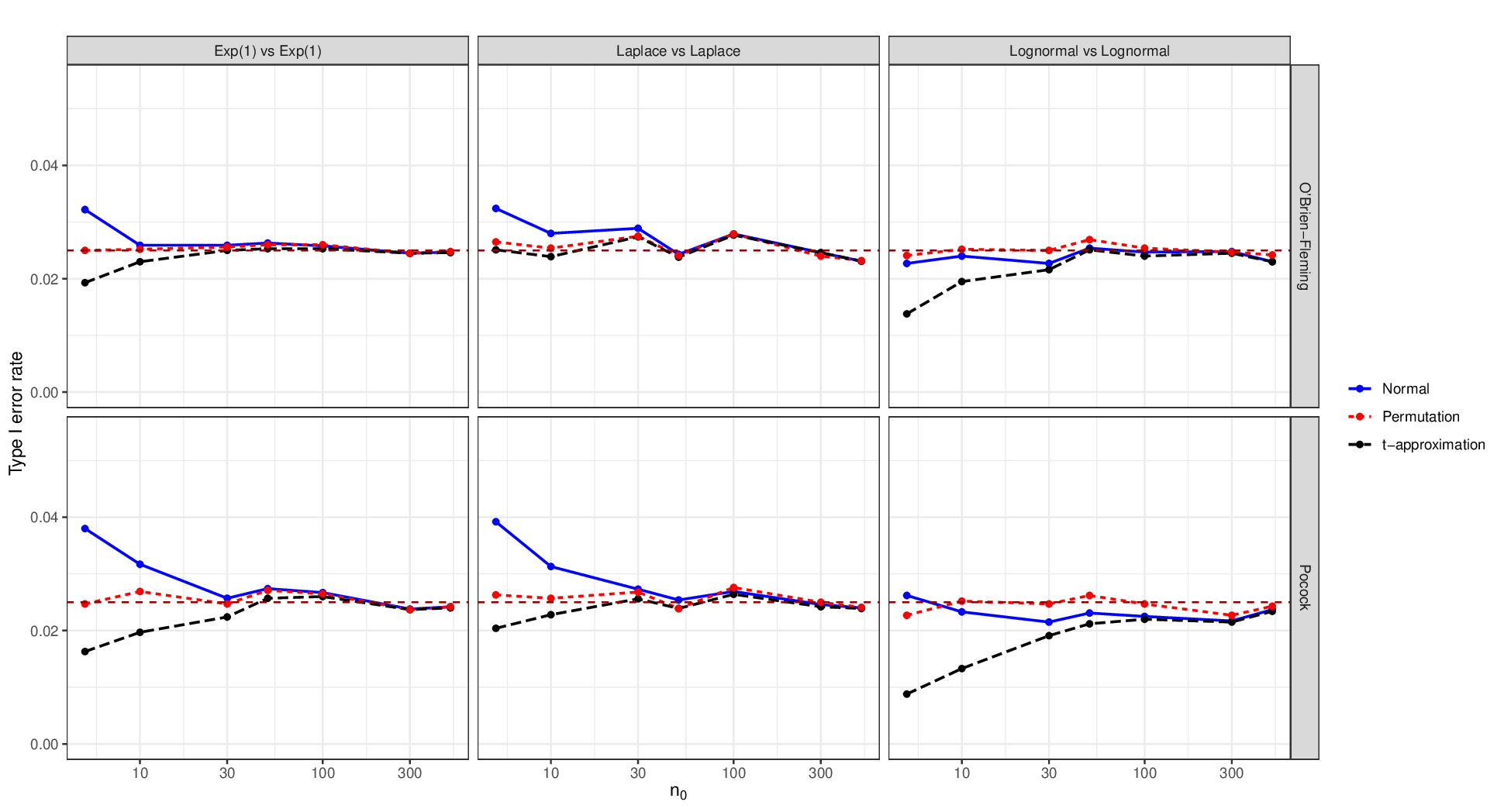}\\
  \caption{Type I error rate ($\alpha=0.025$) simulation results ($y$-axis) of the classical normal method, its $t$-approximation and the proposed permutation method for various scenarios, sample size increments $n_0\in\{5,10,30,50,100,300,500\}$ ($x$-axis), one-sided hypothesis testing, 1:1 allocation ratio and three-stage design}
  \label{three-stages-2}
\end{figure}

\begin{figure}[hpt!]
  \centering
  \includegraphics[width=420pt]{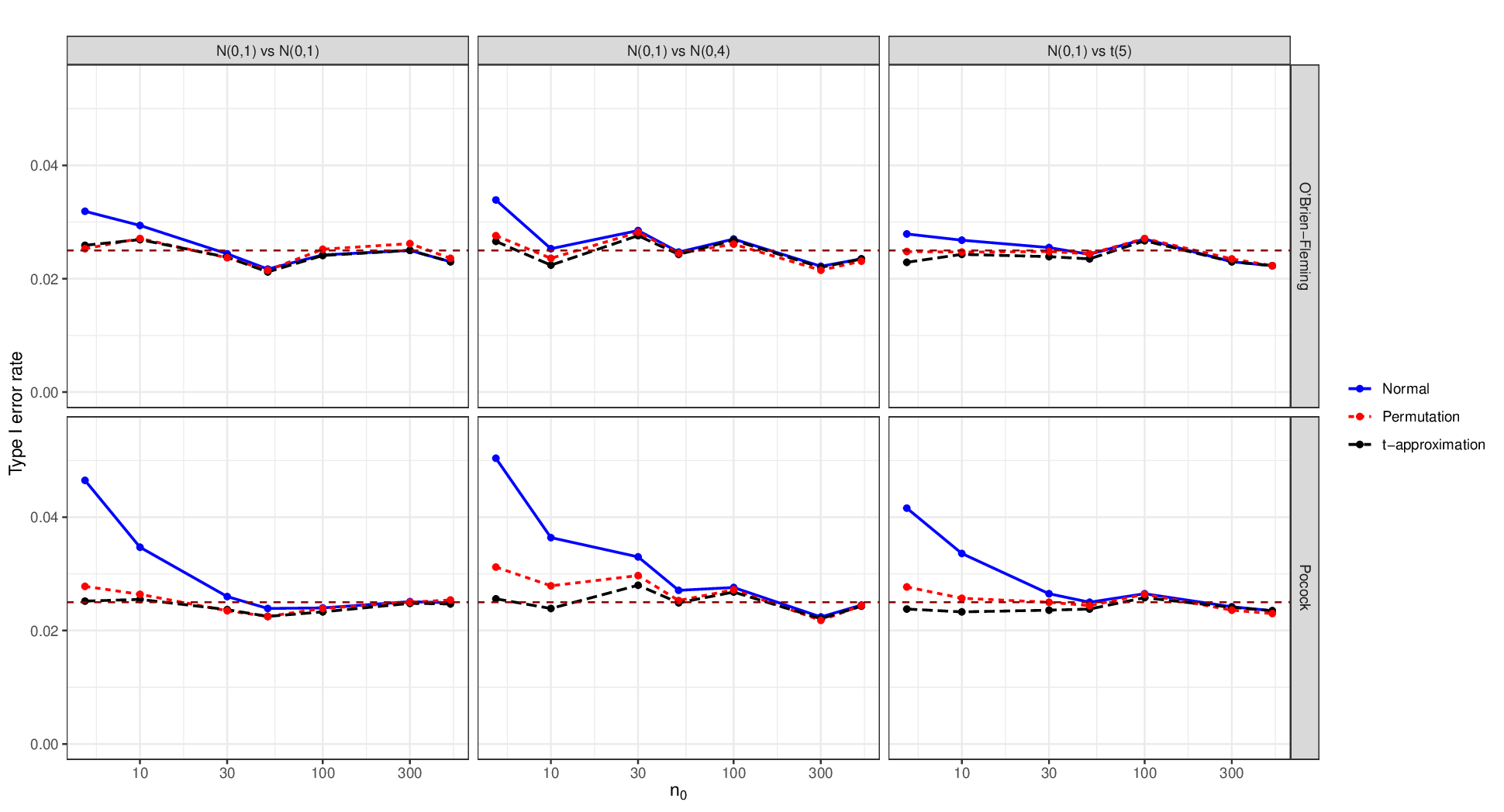}\\
  \caption{Type I error rate ($\alpha=0.025$) simulation results ($y$-axis) of the classical normal method, its $t$-approximation and the proposed permutation method for various scenarios, sample size increments $n_0\in\{5,10,30,50,100,300,500\}$ ($x$-axis), one-sided hypothesis testing, 1:1 allocation ratio and five-stage design}
  \label{five-stages-1}
\end{figure}

\begin{figure}[hpt!]
  \centering
  \includegraphics[width=420pt]{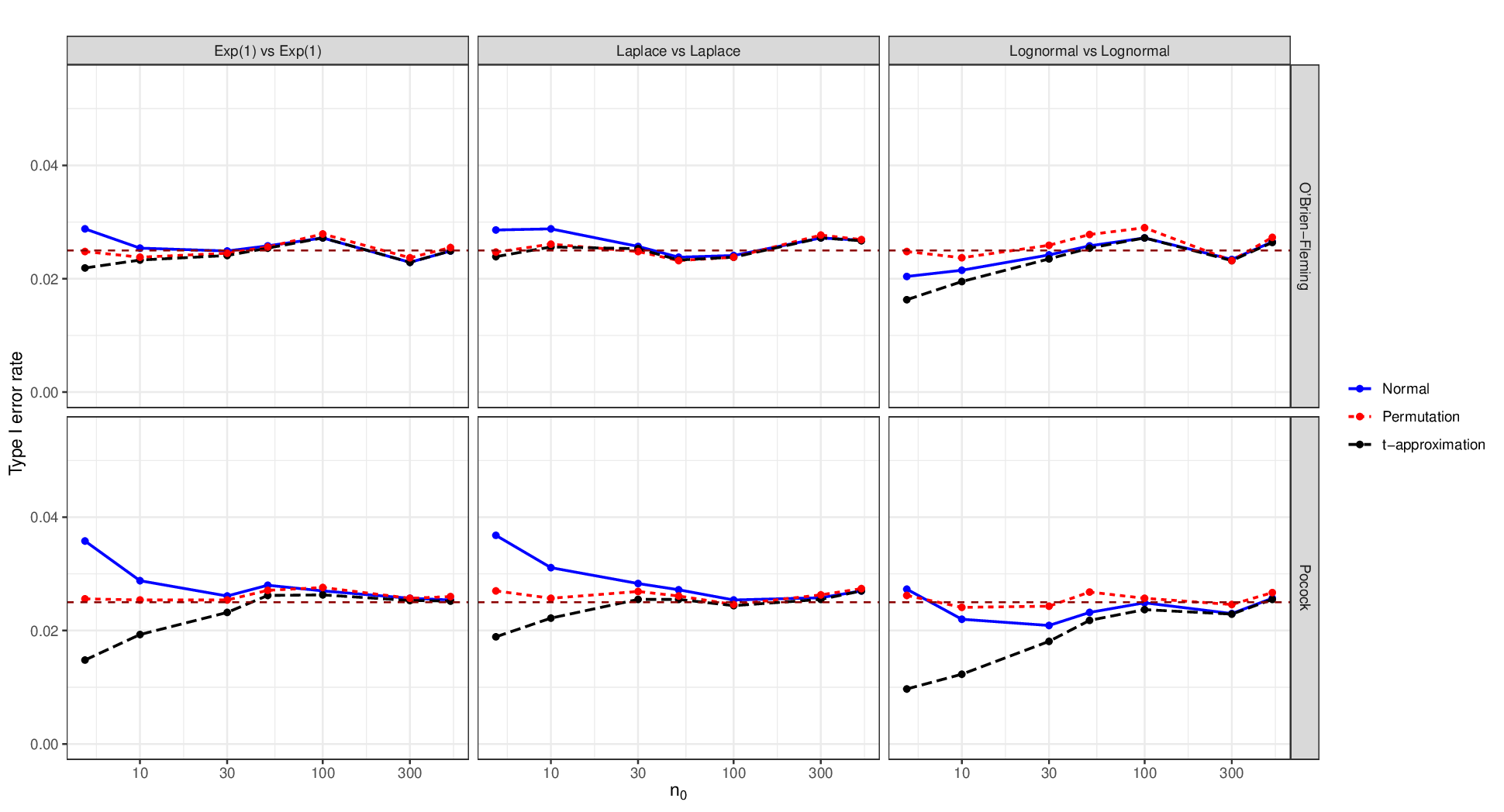}\\
  \caption{Type I error rate ($\alpha=0.025$) simulation results ($y$-axis) of the classical normal method, its $t$-approximation and the proposed permutation method for various scenarios, sample size increments $n_0\in\{5,10,30,50,100,300,500\}$ ($x$-axis), one-sided hypothesis testing, 1:1 allocation ratio and five-stage design}
  \label{five-stages-2}
\end{figure}

The simulation results of the Type I error rate in three-stage and five-stage balanced designs, which are commonly used in medical research, are displayed in Figure \ref{three-stages-1} -- Figure \ref{five-stages-2}.

These results show that the classical group sequential design and its $t$-approximation fail to control the Type I error rate in small sample sizes. With large sample sizes, both the classical method and its $t$-approximation are adequate for the simulated designs. We notice that the $t$-approximation method outperforms the classical method when the data follows a normal or $t$-distribution. However, for skewed distributions, such as the exponential and log-normal distributions, the $t$-approximation method tends to yield conservative results.

The studentized permutation method controls the Type I error rate very well and greatly improves the existing methods. It should be noted that the studentized permutation test serves as a nearly $\alpha$-level test even in the case of small samples for the balanced designs.

\subsection{Power}

\begin{figure}[hpt!]
  \centering
  \includegraphics[width=420pt]{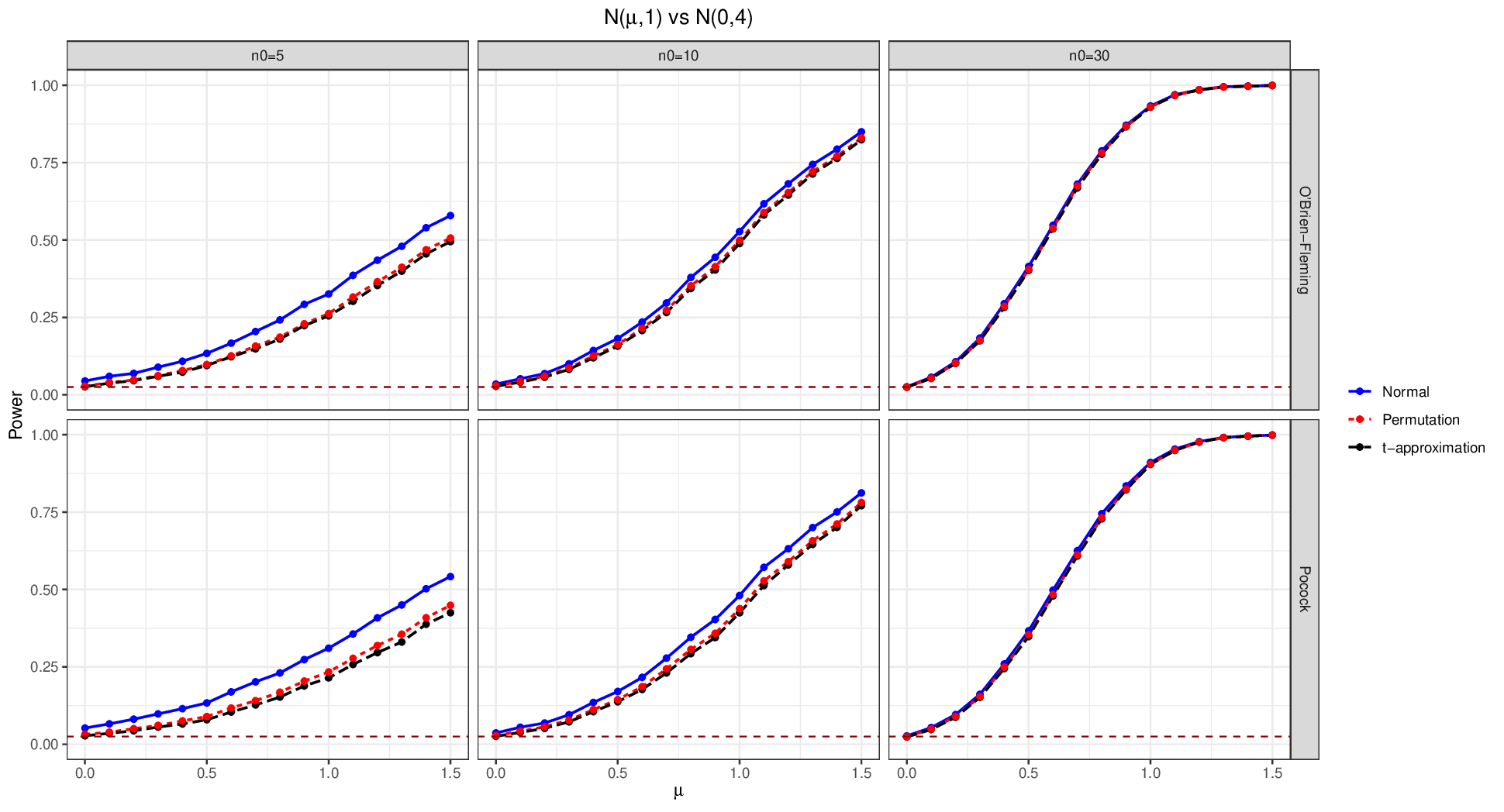}\\
  \caption{Power simulation results ($y$-axis) of the classical normal method, its $t$-approximation and the proposed permutation method for the scenario $N(\mu,1)$ vs $N(0,4)$ under the alternative $\mu\in\{0,0.1,0.2,\ldots,1.5\}$ ($x$-axis), sample size $n_0=5,10,30$, one-sided hypothesis testing, 1:1 allocation ratio and two-stage design}
  \label{power-2}
\end{figure}

\begin{figure}[hpt!]
  \centering
  \includegraphics[width=420pt]{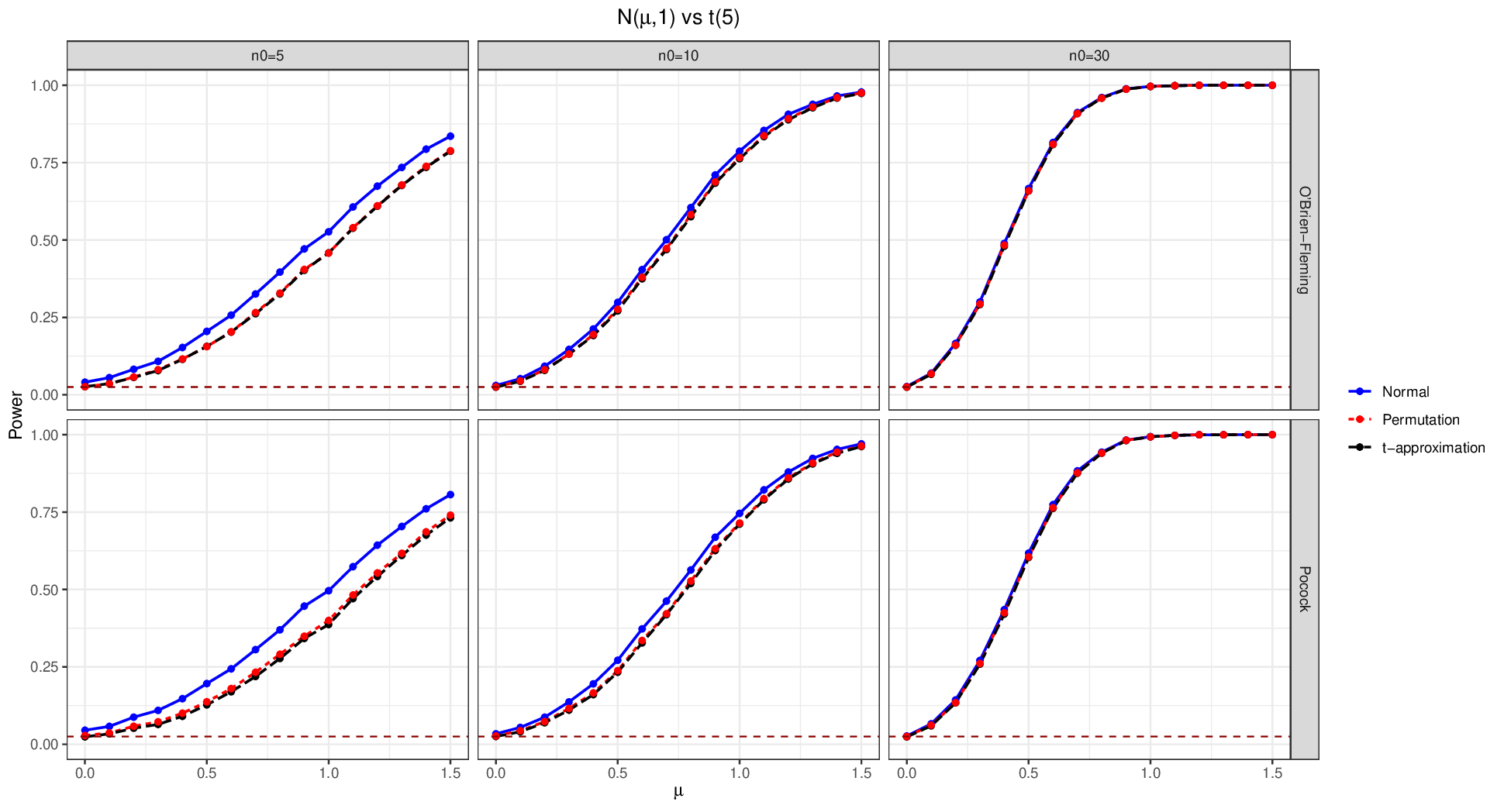}\\
  \caption{Power simulation results ($y$-axis) of the classical normal method, its $t$-approximation and the proposed permutation method for the scenario $N(\mu,1)$ vs $t(5)$ under the alternative $\mu\in\{0,0.1,0.2,\ldots,1.5\}$ ($x$-axis), sample size $n_0=5,10,30$, one-sided hypothesis testing, 1:1 allocation ratio and two-stage design}
  \label{power-3}
\end{figure}

\begin{figure}[hpt!]
  \centering
  \includegraphics[width=420pt]{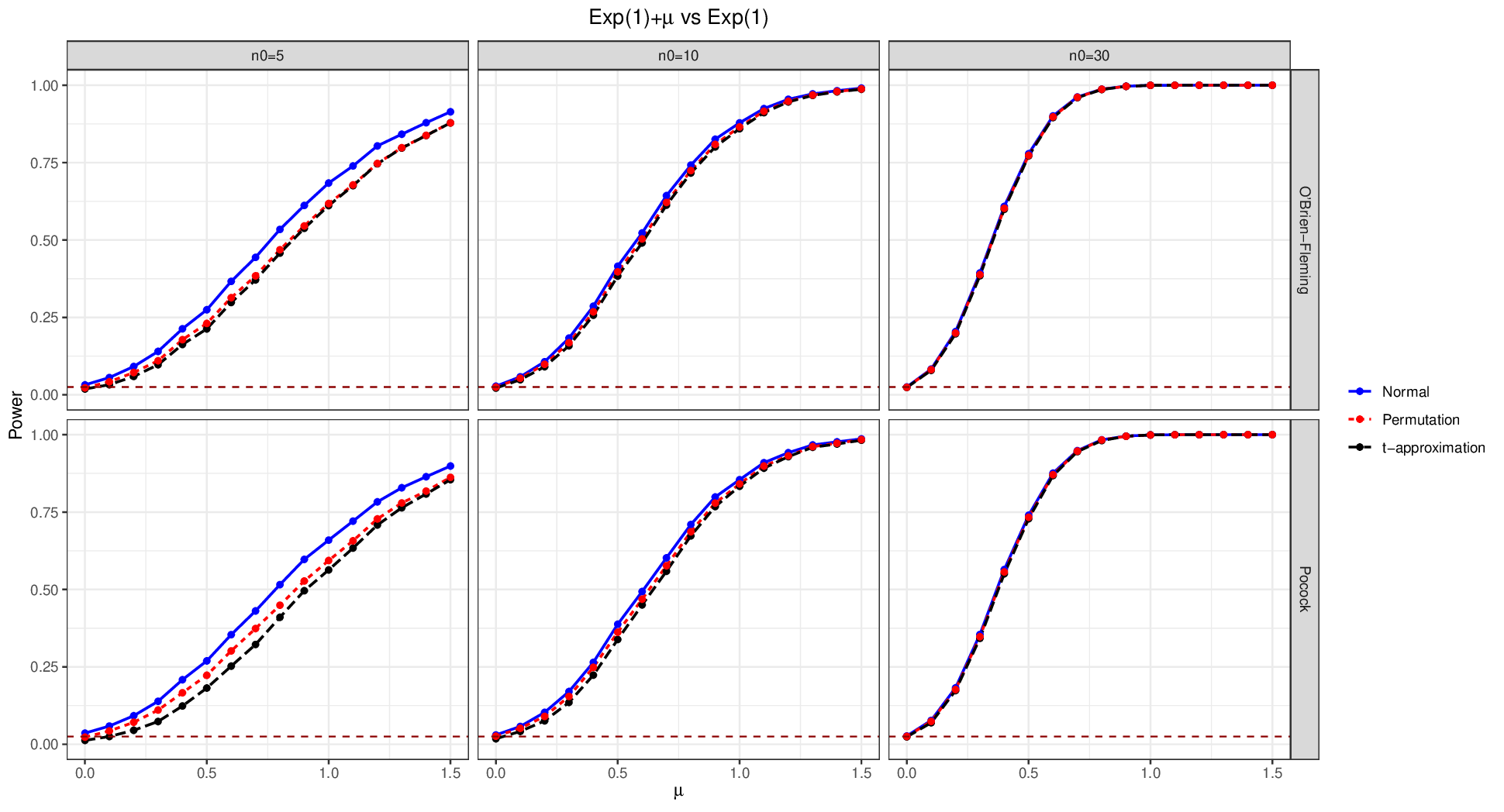}\\
  \caption{{Power simulation results ($y$-axis) of the classical normal method, its $t$-approximation and the proposed permutation method for the scenario Exp(1)+$\mu$ vs Exp(1) under the alternative $\mu\in\{0,0.1,0.2,\ldots,1.5\}$ ($x$-axis), sample size $n_0=5,10,30$, one-sided hypothesis testing, 1:1 allocation ratio and two-stage design}}
  \label{power-4}
\end{figure}

\begin{figure}[hpt!]
  \centering
  \includegraphics[width=420pt]{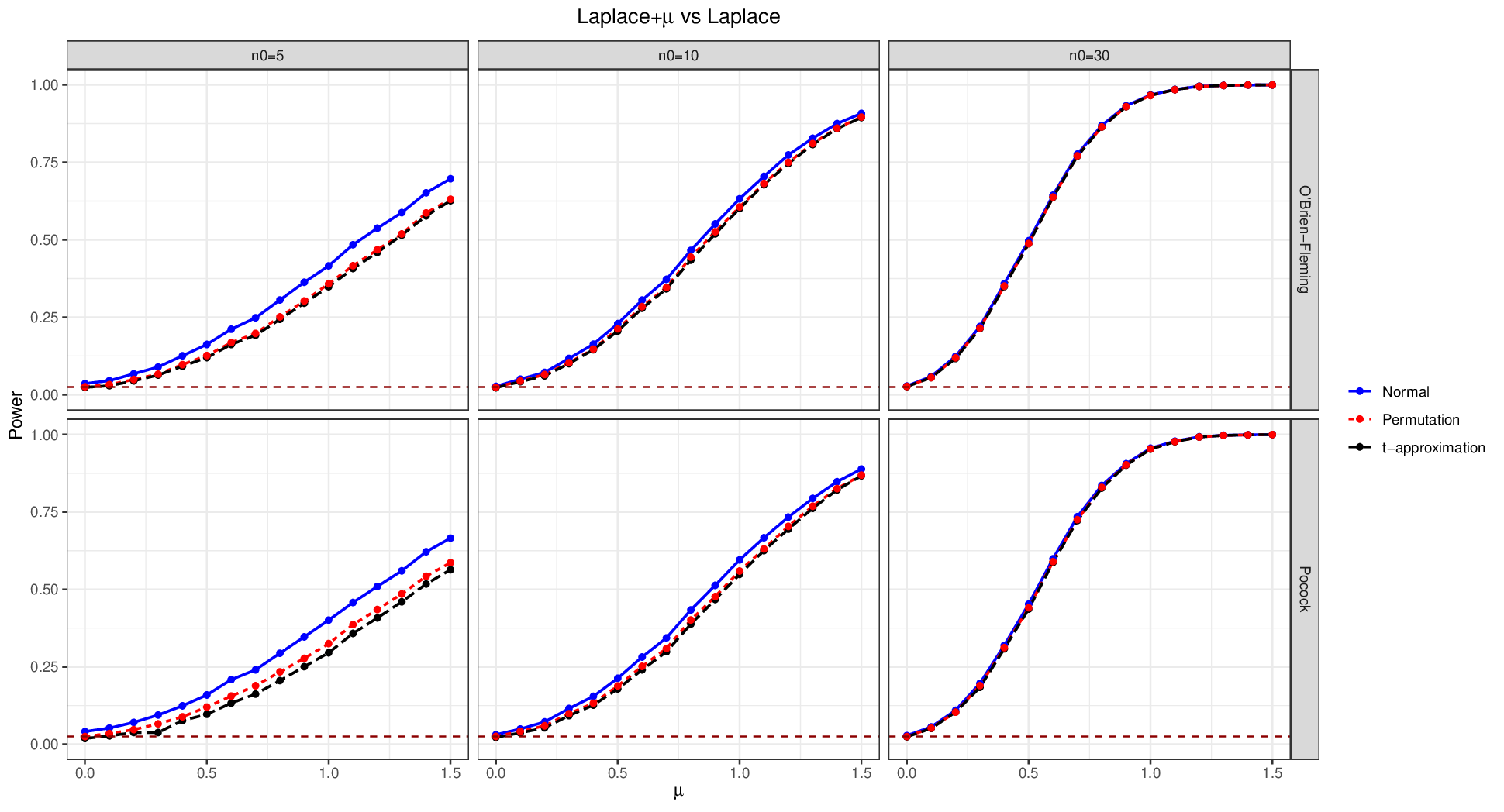}\\
  \caption{{Power simulation results ($y$-axis) of the classical normal method, its $t$-approximation and the proposed permutation method for the scenario Laplace+$\mu$ vs Laplace under the alternative $\mu\in\{0,0.1,0.2,\ldots,1.5\}$ ($x$-axis), sample size $n_0=5,10,30$, one-sided hypothesis testing, 1:1 allocation ratio and two-stage design}}
  \label{power-5}
\end{figure}

\begin{figure}[hpt!]
  \centering
  \includegraphics[width=420pt]{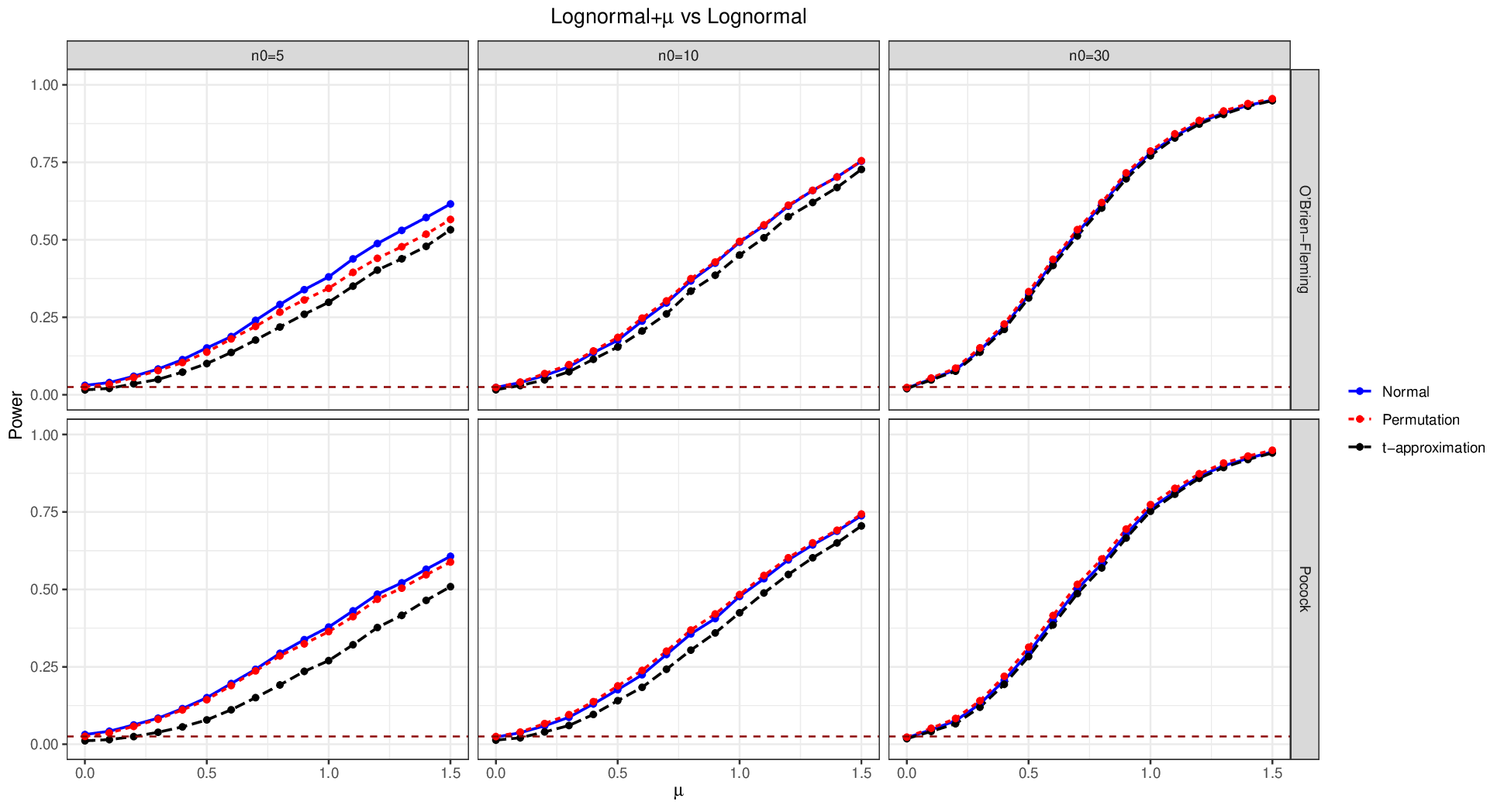}\\
  \caption{{{Power simulation results ($y$-axis) of the classical normal method, its $t$-approximation and the proposed permutation method for the scenario Lognormal+$\mu$ vs Lognormal under the alternative $\mu\in\{0,0.1,0.2,\ldots,1.5\}$ ($x$-axis), sample size $n_0=5,10,30$, one-sided hypothesis testing, 1:1 allocation ratio and two-stage design}}}
  \label{power-6}
\end{figure}

To investigate the power of the procedures, we give a shift $\mu$ for the treatment arms and keep the control arm fixed. To illustrate the power in the case of small sample size, we select $n_0=5,10,30$, where $n_0$ is the sample size in the control arm at any given stage. The simulation results are displayed in Figure \ref{power-2} -- Figure \ref{power-6}.

The power curves from the three competing procedures coincide asymptotically, i.e. no differences among the procedures can be detected with large sample sizes.
The classical test appears to have the best power behaviour owing to its extremely liberal behaviour. Here its power function is not really comparable and we have only included this classical test in the power simulation study for illustration and completeness.

\end{document}